% ceci est la version finale de l'exposee 989
% Titre: "The Strong $ABC$ conjecture over function fields (after McQuillan and Yamanoi)
% 3 juin 2008

\let\noarrow = t
\input eplain

\let\noarrow = t

\input eplain

%%%%%      Magnification, Page dimensions, Margins, Vertical spacing

\magnification=\magstep1

\topskip1truecm
\def\pagewidth#1{
  \hsize=#1
}

\def\pageheight#1{
  \vsize=#1
}

\pageheight{23.5truecm} \pagewidth{16truecm}

\abovedisplayskip=3mm \belowdisplayskip=3mm
\abovedisplayshortskip=0mm \belowdisplayshortskip=2mm
\parindent1pc

\normalbaselineskip=13pt \baselineskip=13pt

\def\spacing{{\smallskip}}

%%%%%%%%%%%%%%   The following moves the actual dvi file
%%%%%%%%%%%%%%  \voffset vertically and \hoffset horizontally
\voffset=0pc \hoffset=0pc

%%%%%%%%%%%%%%%  To be used in  the abstract

\newdimen\abstractmargin
\abstractmargin=3pc

%%%%%%%%%%%%%  To be used in footnote first page

\newdimen\footnotemargin
\footnotemargin=1pc

%%%%%%%%%%%%  Fonts for the abstract from Knuth's manual

\font\eightrm=cmr8 \relax % preloaded in plain.tex
\font\sixrm=cmr6 \relax % preloaded in plain.tex
\font\eighti=cmmi8 \relax     \skewchar\eighti='177 % preloaded
\font\sixi=cmmi6 \relax       \skewchar\sixi='177   % preloaded
\font\eightsy=cmsy8 \relax    \skewchar\eightsy='60 % preloaded
\font\sixsy=cmsy6 \relax      \skewchar\sixsy='60   % preloaded
\font\eightbf=cmbx8 \relax % preloaded in plain.tex
\font\sixbf=cmbx6 \relax   % preloaded in plain.tex
\font\eightit=cmti8 \relax % preloaded in plain.tex
\font\eightsl=cmsl8 \relax % preloaded in plain.tex
\font\eighttt=cmtt8 \relax % preloaded in plain.tex

\catcode`\@=11
\newskip\ttglue

\def\eightpoint{\def\rm{\fam0\eightrm}%
 \textfont0=\eightrm \scriptfont0=\sixrm
 \scriptscriptfont0=\fiverm
 \textfont1=\eighti \scriptfont1=\sixi
 \scriptscriptfont0=\fivei
 \textfont2=\eightsy \scriptfont2=\sixsy
 \scriptscriptfont2=\fivesy
 \textfont3=\tenex \scriptfont3=\tenex
 \scriptscriptfont3=\tenex
 \textfont\itfam\eightit \def\it{\fam\itfam\eightit}%
 \textfont\slfam\eightsl \def\sl{\fam\slfam\eightsl}%
 \textfont\ttfam\eighttt \def\tt{\fam\ttfam\eighttt}%
 \textfont\bffam\eightbf \scriptfont\bffam\sixbf
   \scriptscriptfont\bffam\fivebf \def\bf{\fam\bffam\eightit}%
 \tt \ttglue=.5em plus.25em minus.15em
 \normalbaselineskip=9pt
 \setbox\strutbox\hbox{\vrule height7pt depth3pt width0pt}%
 \let\sc=\sixrm \let\big=\eifgtbig \normalbaselines\rm}

%%%%%%%              Fonts

 \font\titlefont=cmbx12 scaled\magstep1
 \font\sectionfont=cmbx12
 \font\ssectionfont=cmsl10
 \font\claimfont=cmsl10

 \font\normalfont=cmr10

%%     *****Part of MSSYMB.Tex*****
%%%%%      Eulerscript
\catcode`\@=11 \font\teneusm=eusm10 
\font\seveneusm=eusm7  \font\fiveeusm=eusm5
\newfam\eusmfam \textfont\eusmfam=\teneusm
\scriptfont\eusmfam=\seveneusm \scriptscriptfont\eusmfam=\fiveeusm
\def\hexnumber@#1{\ifcase#1
0\or1\or2\or3\or4\or5\or6\or7\or8\or9\or         A\or B\or C\or D\or
E\or F\fi } \edef\eusm@{\hexnumber@\eusmfam}
\def\euscr{\ifmmode\let\next\euscr@\else
\def\next{\errmessage{Use \string\euscr\space only in math mode}}\fi\next}
\def\euscr@#1{{\euscr@@{#1}}} \def\euscr@@#1{\fam\eusmfam#1} \catcode`\@=12

%%     *****Part of MSSYMB.Tex*****
%%%%%      Euler
\catcode`\@=11 \font\teneuex=euex10 
 \font\seveneuex=euex7  \newfam\euexfam
\textfont\euexfam=\teneuex  \scriptfont\euexfam=\seveneuex
 \def\hexnumber@#1{\ifcase#1
0\or1\or2\or3\or4\or5\or6\or7\or8\or9\or         A\or B\or C\or D\or
E\or F\fi } \edef\euex@{\hexnumber@\euexfam}
\def\euscrex{\ifmmode\let\next\euscrex@\else
\def\next{\errmessage{Use \string\euscrex\space only in math mode}}\fi\next}
\def\euscrex@#1{{\euscrex@@{#1}}} \def\euscrex@@#1{\fam\euexfam#1}
\catcode`\@=12

%%     *****Part of MSSYMB.Tex*****
%%%%%      Eulerscript
\catcode`\@=11 \font\teneufb=eufb10 
\font\seveneufb=eufb7  \font\fiveeufb=eufb5
\newfam\eufbfam \textfont\eufbfam=\teneufb
\scriptfont\eufbfam=\seveneufb \scriptscriptfont\eufbfam=\fiveeufb
\def\hexnumber@#1{\ifcase#1
0\or1\or2\or3\or4\or5\or6\or7\or8\or9\or         A\or B\or C\or D\or
E\or F\fi } \edef\eufb@{\hexnumber@\eufbfam}
\def\euscrfb{\ifmmode\let\next\euscrfb@\else
\def\next{\errmessage{Use \string\euscrfb\space only in math mode}}\fi\next}
\def\euscrfb@#1{{\euscrfb@@{#1}}} \def\euscrfb@@#1{\fam\eufbfam#1}
\catcode`\@=12

%%     *****Part of MSSYMB.Tex*****
%%%%%      Eulerscript
\catcode`\@=11 \font\teneufm=eufm10 
\font\seveneufm=eufm7  \font\fiveeufm=eufm5
\newfam\eufmfam \textfont\eufmfam=\teneufm
\scriptfont\eufmfam=\seveneufm \scriptscriptfont\eufmfam=\fiveeufm
\def\hexnumber@#1{\ifcase#1
0\or1\or2\or3\or4\or5\or6\or7\or8\or9\or         A\or B\or C\or D\or
E\or F\fi } \edef\eufm@{\hexnumber@\eufmfam}
\def\euscrfm{\ifmmode\let\next\euscrfm@\else
\def\next{\errmessage{Use \string\euscrfm\space only in math mode}}\fi\next}
\def\euscrfm@#1{{\euscrfm@@{#1}}} \def\euscrfm@@#1{\fam\eufmfam#1}
\catcode`\@=12

%%     *****Part of MSSYMB.Tex*****
%%%%%      Eulerscript
\catcode`\@=11 \font\teneusb=eusb10 
\font\seveneusb=eusb7  \font\fiveeusb=eusb5
\newfam\eusbfam \textfont\eusbfam=\teneusb
\scriptfont\eusbfam=\seveneusb \scriptscriptfont\eusbfam=\fiveeusb
\def\hexnumber@#1{\ifcase#1
0\or1\or2\or3\or4\or5\or6\or7\or8\or9\or         A\or B\or C\or D\or
E\or F\fi } \edef\eusb@{\hexnumber@\eusbfam}
\def\euscrsb{\ifmmode\let\next\euscrsb@\else
\def\next{\errmessage{Use \string\euscrsb\space only in math mode}}\fi\next}
\def\euscrsb@#1{{\euscrsb@@{#1}}} \def\euscrsb@@#1{\fam\eusbfam#1}
\catcode`\@=12

%%%%%%%%  \Bbb
%%
\catcode`\@=11 \font\tenmsa=msam10 
\font\sevenmsa=msam7  \font\fivemsa=msam5
\font\tenmsb=msbm10  \font\sevenmsb=msbm7
 \font\fivemsb=msbm5 \newfam\msafam
\newfam\msbfam \textfont\msafam=\tenmsa
\scriptfont\msafam=\sevenmsa
  \scriptscriptfont\msafam=\fivemsa
\textfont\msbfam=\tenmsb  \scriptfont\msbfam=\sevenmsb
  \scriptscriptfont\msbfam=\fivemsb
\def\hexnumber@#1{\ifcase#1 0\or1\or2\or3\or4\or5\or6\or7\or8\or9\or
        A\or B\or C\or D\or E\or F\fi }
\edef\msa@{\hexnumber@\msafam} \edef\msb@{\hexnumber@\msbfam}
\mathchardef\square="0\msa@03 \mathchardef\subsetneq="3\msb@28
\mathchardef\supsetneq="3\msb@29 \mathchardef\ltimes="2\msb@6E
\mathchardef\rtimes="2\msb@6F \mathchardef\dabar="0\msa@39
\mathchardef\daright="0\msa@4B \mathchardef\daleft="0\msa@4C
\def\dashrightarrow{\mathrel{\dabar\dabar\daright}}

\def\Bbb{\ifmmode\let\next\Bbb@\else
        \def\next{\errmessage{Use \string\Bbb\space only in math mode}}\fi\next}
\def\Bbb@#1{{\Bbb@@{#1}}}
\def\Bbb@@#1{\fam\msbfam#1}
\catcode`\@=12

%%
%%      *****End of MSSYMB.Tex*****

%%%%%%%           Environments

%%% Section,Subsection,Paragraph numbers
\newcount\senu
\def\senum{\number\senu}
\newcount\ssnu
\def\ssnum{\number\ssnu}
\newcount\fonu
\def\fonum{\number\fonu}

\def\num{{\senum.\ssnum}}
\def\numfo{{\senum.\ssnum.\fonum}}

%%% Section, Subsection, Bibliography enviroment

\outer\def\section#1\par{\vskip0pt
  plus.3\vsize\penalty20\vskip0pt
  plus-.3\vsize\bigskip\vskip\parskip
  \message{#1}\centerline{\sectionfont\senum\enspace#1.}
  \nobreak\smallskip}

\def\endsection{\advance\senu by1\penalty-20\smallskip\ssnu=1}
\outer\def\ssection#1\par{\bigskip
  \message{#1}{\noindent\bf\num\ssectionfont\enspace#1.\thinspace}
  \nobreak\normalfont}

\def\endssection{\advance\ssnu by1\smallskip\ifdim\lastskip<\medskipamount
\removelastskip\penalty55\medskip\fi\fonu=1\normalfont}

%%%%%% Theorem, Corollary, Lemma, Remark, Proposition, Definition
\def\proclaim #1\par{\bigskip
  \message{#1}{\noindent\bf\num\enspace#1.\thinspace}
  \nobreak\claimfont}

\def\defi{\proclaim Definition\par}
\def\lemma{\proclaim Lemma\par}
\def\prop{\proclaim Proposition\par}
\def\rmk{\proclaim Remark\par\normalfont}
\def\thm{\proclaim Theorem\par}
\def\conj{\proclaim Conjecture\par}

\def\enddefi{\endssection}
\def\endlemma{\endssection}
\def\endprop{\endssection}
\def\endrmk{\endssection}
\def\endthm{\endssection}
\def\endconj{\endssection}

\def\Proof{{\noindent\sl Proof: \/}}

%%%%%%%%%%%       Arrows

\def\maplefto#1{\ \smash{\mathop{\longleftarrow}\limits^{#1}}\ }

\def\llongrightarrow{\relbar\joinrel\relbar\joinrel\rightarrow}
\def\lllongrightarrow{\hbox to 40pt{\rightarrowfill}}

\def\twoheadrightarrow{\rightarrow\kern -8pt\rightarrow}

\def\maprighto#1{\smash{\mathop{\longrightarrow}\limits^{#1}}}

\def\mapdownr#1{\Big\downarrow\rlap{$\vcenter{\hbox{$\scriptstyle#1$}}$}}
\def\mapdownl#1{\llap{$\vcenter{\hbox{$\scriptstyle#1$}}$}\Big\downarrow}

\def\llongmaprighto#1{\ \smash{\mathop{\llongrightarrow}\limits^{#1}}\ }

\def\lllongmaprighto#1{\ \smash{\mathop{\lllongrightarrow}\limits^{#1}}\ }

\def\longleftmapsto{\longleftarrow\kern-2pt\mapstochar\;}%%

\def\llongmapsto{\,\vert\kern-3.2pt\joinrel\longrightarrow\,}
\def\llongmapsto{\,\vert\kern-3.7pt\joinrel\llongrightarrow\,}
\def\lllongmapsto{\,\vert\kern-5.5pt\joinrel\lllongrightarrow\,}

\def\isomarrow{\maprighto{\lower3pt\hbox{$\scriptstyle\sim$}}}
\def\llongisomarrow{\llongmaprighto{\lower3pt\hbox{$\scriptstyle\sim$}}}
\def\lllongisomarrow{\lllongmaprighto{\lower3pt\hbox{$\scriptstyle\sim$}}}

\def\lisomarrow{\maplefto{\lower3pt\hbox{$\scriptstyle\sim$}}}

%%%%%%%%%%        References
%%%  Prints on the margin the label.
\font\labprffont=cmtt8
\def\strutdepth{\dp\strutbox}
\def\labtekst#1{\vtop to \strutdepth{\baselineskip\strutdepth\vss\llap{{\labprffont #1}}\null}}
\def\marglabel#1{\strut\vadjust{\kern-\strutdepth\labtekst{#1\ }}}

\def\label #1. #2\par{{\definexref{#1}{\num}{#2}}}
\def\labelf #1\par{{\definexref{#1}{\numfo}{formula}}}
\def\labelse #1\par{{\definexref{#1}{\num}{section}}}

%%%%%%%%%%%%%  Math relations

\def\*{{\ast}}

\def\QQ{{\Bbb Q}}
\def\ZZ{{\Bbb Z}}
\def\RR{{\Bbb R}}

\def\CC{{\Bbb C}}
\def\PP{{\Bbb P}}
\def\LL{{\Bbb L}}
\def\MM{{\Bbb M}}
\def\cO{{\cal O}}
\def\II{{\Bbb I}}

\def\Proj{{\rm Proj}}
\def\Spec{{\rm Spec}}

\def\p{{\euscrfm p}}

\def\inf{{\rm inf}}

\def\SGA1{{\rm SGA1}}

\def\M{{\cal M}}
\def\U{{\cal U}}

\senu=1 \ssnu=1 \fonu=1

\noindent{S\'eminaire BOURBAKI}

\noindent{60\'{e}me ann\'ee, 2007--2008, ${\rm n}^\circ$ 989}

\

\centerline{\titlefont The strong $ABC$ conjecture over function
fields}

\centerline{\bf [after McQuillan and Yamanoi]}

\

\spacing \centerline{\bf C. Gasbarri}

\bigskip

\section Introduction\par

\

One of the deepest conjecture in arithmetic is the $abc$ conjecture:

\label abc. conjecture\par\conj Let $\epsilon >0$, then there exists
a constant $C(\epsilon)$ for which the following holds: Let $a$, $b$
and $c$ three integral numbers such that $(a,b)=1$ and $a+b=c$ then
$$\max\{ \vert a\vert, \vert b\vert, \vert c\vert\}\leq
C(\epsilon)\left(\prod_{p/abc}p\right)^{1+\epsilon}$$ where the
product is taken over all the prime numbers dividing $abc$.
\endconj

Let's give a geometric interpretation of this conjecture:

Consider the arithmetic surface $\PP^1_{\ZZ}\to\Spec(\ZZ)$ equipped
with the tautological line bundle $\cO(1)$ and the divisor
$D:=[0:1]+[1:0]+[1:-1]$. Suppose we have a section $P:\Spec(\ZZ)\to
\PP^1_{\ZZ}$, not contained in $D$, then $P^\*(D)$ is an effective
Weil divisor on $\Spec(\ZZ)$ which can be written as
$\sum_pv_p(D)[p]$.

Define {\it the radical} of the divisor as
$N_D^{(1)}(P):=\sum_p\min(1,v_p(D))\log(p)$.

The conjecture can be stated in this way: for every $\epsilon>0$
there is a constant $C(\epsilon)$ such that, for every section
$P:\Spec(\ZZ)\to\PP^1_{\ZZ}$ we have
$$h_{\cO(1)}(P)\leq (1+\epsilon)N^{(1)}_D(P)+C(\epsilon).$$
Where $h_{\cO(1)}(P)$ is the height of $P$ with respect to
$\cO_{\PP^1}(1)$. When we state the conjecture in this way we see
many possible generalizations. We also clearly see the geometric
analogue over function fields (cf. next sections for details). Let's
formulate the conjecture in the most general version.

If $K$ is a number field, we denote by $\cO_K$ the ring of integers
of $K$ and $\Delta_K$ its discriminant. If $X\to K$ is an arithmetic
surface. $D$ is an effective divisor over $X$ and $P:\Spec(\cO_K)\to
X$, not contained in $D$, we define the {\it radical} of $D$ as the
real number
$N_D^{(1)}(P):=\sum_{\p\in\Spec\max(\cO_K)}\min\{1;v_{\p}P^\*(D)\}\log
Card(\cO_K/\p)$. The  general strong $abc$ conjecture is the
following:

\label generalabc. conjecture\par\conj Let $\epsilon>0$, and $K$ be
a number field, $\pi:X\to\Spec(\cO_K)$ a regular arithmetic surface
and $D\hookrightarrow X$ an effective divisor on $X$. Denote by
$K_{X/\cO_K}$ the relative dualizing sheaf. Then there exists a
constant $C:=C(X,\epsilon, D)$ for which the following holds: let
$L$ be a finite extension of $K$ and $P:\Spec(\cO_L)\to X$, not
contained in $D$, then
$$h_{K_{X/\cO_K}(D)}(P)\leq
(1+\epsilon)(N^{(1)}_D(P)+\log\vert\Delta_L\vert)+C[L:K]$$

where $h_{K_{X/\cO_K}(D)}(P)$  is the height of $P$ with respect to
$K_{X/\cO_K}(D)$.
\endconj

We will not list here the endless number of consequences of this
conjecture and we refer to [BG] or to the web page [NI] for details.
One may also see the report [OE] in this seminar. We only notice
that, if such a conjecture was true, more or less all the possible
problems about the arithmetic of algebraic curves over number fields
would have an effective answer: for instance one easily sees that, if the constant $C$ is effective, it
easily implies the famous Fermat Last Conjecture (now a theorem
[WI]) and it allows to solve effectively diophantine equations in two variables:

\label diophantine. theorem\par\thm Suppose that conjecture \ref{generalabc} is true. Let $F(x;y)\in\ZZ[x,y]$ be an irreducible polynomial of degree at least three. Then there exists a constant $C$, depending only on $F$, such that for every number field $K$ and for every solution $(x;y)\in\cO_K\times\cO_K$ of the diophantine equation
$F(x;y)=0$, we have
$$h_{\cO(1)}([x:y:1])\leq (1+\epsilon)\log\vert\Delta_K\vert +C_\epsilon[K:\QQ].$$
In particular there are only finitely many  solutions in $\cO_K\times\cO_K$ and
their height can be explicitly bounded.
\endthm

Observe that, if the conjecture is true and the constant
$C_\epsilon$ is explicit, then we can explicitly {\it compute and
find} the set of solutions of the diophantine equation in $\cO_K\times\cO_K$.

Similarly we may obtain an effective version of Mordell conjecture (Faltings theorem) and of the classical Siegel
theorem on integral points of hyperbolic curves.

At the moment we know that the set of integral points of an
hyperbolic curve  is finite (projective by Faltings theorem [FA] or
affine by Siegel theorem cf. [SE]) but we are not able to explicitly
bound their height (up to some sporadic cases); thus, in particular,
it is not possible to find all the rational points of an hyperbolic
curve.

In this paper we will report about the solution of the analogue of
the $abc$ conjecture over function fields (for the analogy between
number fields and function fields arithmetic cf. for instance [SE]).

The analogue of conjecture \ref{abc} for polynomials is quite easy
and proved in [MA]: If $f$ is a polynomial over $\CC$ (to simplify),
let  $N_0=(f)$ be the number of {\it distinct} roots of $f$. Then
the analogue of the $abc$ conjecture for polynomials is

\label abcforpoly. theorem\par\thm (Mason) Let $f$, $g$ and $h$
three polynomials relatively coprime in $\CC[t]$ such that $f+g=h$,
then
$$\max\{\deg(f), \deg(g),\deg(h)\}\leq N_0(fgh)-1.$$
\endthm

This theorem is the analogue of conjecture \ref{generalabc} for
function fields when $X=\PP^1\times\PP^1$, $\pi:X\to\PP^1$ is the
first projection,
$D=\PP^1\times[0:1]+\PP^1\times[1:0]+\PP^1\times[1:-1]$ and $P$ is a
section.  One easily deduce it from Hurwitz formula (cf. next
section). It can be seen as the beginning of all the story, and it
has some interesting consequences: for example the analogue of
Fermat last theorem for polynomials is an immediate consequence of
it. Usually statements in the function fields situation are much
easier to prove then their correspondent in the number fields
situation. In this case one should notice an amazing point: Suppose
that, over number fields, we can prove conjecture \ref{generalabc}
when $X=\PP^1_\ZZ$
 and $D=[0:1]+[1:0]+[1:-1]$ then we
can deduce from this the general case! To prove this one applies the
proof of theorem \ref{reduction} to a suitable Belyi map (for more
details cf. [EL]). In the function fields case this is {\it not} the
case! We cannot deduce the general case from an isotrivial case. For
this reason it is our opinion that $\PP^1_\ZZ$ with the divisor
$[0:1]+[1:0]+[1:-1]$ (unit equations) is a highly non isotrivial
family over $\Spec(\ZZ)$ (whatever an isotrivial family should be).

Exploiting the analogy between the arithmetic geometry over number
fields and the theory of {\it analytic} maps from a parabolic curve
to a surface (cf. for instance [VO1]), an analogue of the $abc$
conjecture for these maps also is solved.

We will propose two proof of the $abc$ conjecture over function
fields (and for analytic maps). The first is the proof by McQuillan
[MQ3] and the second is by Yamanoi [YA3]. The proof by McQuillan is
synthetically explained  in the  original paper; it make a
systematic use of the theory of integration on algebraic stacks;
although this is very natural in this contest, it needs a very heavy
background (which here is used only in a quite easy situation). Thus
we preferred to propose a self contained proof which uses the
(easier) theory of normal $\QQ$--factorial varieties; the proof
follows the main ideas of the original one. The proof by Yamanoi
requires skillful combinatorial computations, well explained in the
original paper, thus we preferred to sketch his proof in a special
(but non trivial) case: the main ideas and tools are all used and we
think that once one understand this case, it is easier to follow the
proof in the general situation.

As before, as a consequence we find, for instance,  a strong
effective version of Mordell conjecture over function fields (in
characteristic zero), for non isotrivial families of hyperbolic
curves.

In the next section we will explain why the $abc$ conjecture for
isotrivial curves corresponds respectively to the Hurwitz  formula
in the geometric case and to the Nevanlinna Second Main theorem in
the analytic case. Thus the $abc$ conjecture may be seen as a non
isotrivial version of these theorems.

There are at least two strategies to attack the Second Main Theorem
of Nevanlinna theory. The first strategy uses tools from analytic
and differential geometry,  it is strictly related to the algebraic
geometry of the Hurwitz formula and to the existence of particular
singular metrics on suitable line bundles: it has been strongly
generalized to analytic maps between equidimensional varieties by
Griffiths, King and others in the 70's (cf. [GK]). The second
strategy is via Ahlfors theory (cf. [AH]) it is  much related to the
algebraic and combinatorial topology of maps between surfaces; the
version of the SMT one obtain int this way is weaker then the
original one but also more subtle: one sees that one can perturb a
little bit the divisor $D$ without perturbing the statement (cf. \S
8). These two approaches correspond respectively to the two proposed
proofs. The Proof by McQuillan is nearer to the first strategy while
the Yamanoi's is more topological. One should notice that, while the
first proof is predominantly of global nature and the second is
essentially local, both meet the main difficulties in an argument
which is localized around the singular points of the morphism
$p:X\to B$. If the morphism $p$ is relatively smooth, McQuillan's
proof is much simpler. In a hypothetical relatively smooth case,
Yamanoi approach reduces to the Ahlfors theory: you will  observe
that, unless you are in the isotrivial case, in the Yamanoi approach
there is  always  bad reduction.

Both proofs holds for curves over function fields in one variable
over $\CC$ and both heavily use analytic and topological methods,
specific of the complex topology. We should notice that the analogue
of the $abc$ conjecture, as stated before, over a function field
with positive characteristic is false! (cf. [KI]).

\ssection A short overview of the history of the $abc$
conjecture.\par The $abc$ conjecture has a weak and a strong version
(in the arithmetic case they are both unproven and very deep). Over
function fields, the weak $abc$ is easier to prove and it is
strictly related with the theory of elliptic curves (cf. [HS] and
[SZ]). Here we deal with the strong version. The conjecture have
been formulated in the middle 80's by Masser and Oesterl\'e
exploiting the analogy between number fields and function fields and
the version for polynomials proved in [MA]. The general version, as
stated here, have been formulated by Vojta in [VO1] as a consequence
of a series of conjectures for varieties of arbitrary dimension. The particular case of $\PP^1\times\PP^1$ and $D:=[0:1]\times\PP^1+[1:0]\times\PP^1+[1:1]\times\PP^1+\Delta$ ($\Delta$ being the diagonal) was previously proposed by Oesterl\'e.
Some weak versions of the conjecture in the contest of the Value
distribution theory have been proved in the papers [SA] and [OS].

In the paper [VO2] one find a proof of a weak version of the
conjecture (with factor $2+\epsilon$ instead of $1+\epsilon$) in the
algebraic case when $D$ is empty (function field case!); it can be
easily generalized to the case when $D$ is arbitrary. In the recent
paper [MY] we can find an algebraic  proof of a weak version of the
conjecture. On the preprint [CH] one can find another overview of
the proofs.

Recently  one can find generalizations of the theory in the papers
[YA1] and [NWY]. A strong generalization (and many other results),
for families of surfaces, is proved in the forthcoming book [MQ4].
\endssection

\ssection Acknowledgements\par I would like to warmly thank Michael
McQuillan for his continuous help and support. I also thank  J. B.
Bost, A. Chambert Loir and P. Vojta for their comments to the paper.

\endssection

\endsection

\section Notations and overview of Value Distribution Theory\par

\

If $X$ is a set and $U\subseteq X$ is a subset, then we denote by
$\II_U$ the characteristic function of $U$.

Let $X$ be a smooth variety defined over $\CC$ and $L$ a line bundle
$D$ a Cartier divisor over $X$ such that $L=\cO_X(D)$. We Suppose
that $L$ is equipped with a smooth metric.

Over $X$ we have the two operators $d:=\partial +\overline\partial$
and $d^c:={{1}\over{4\pi\sqrt{-1}}}(\partial-\overline\partial)$.

The divisor $D$ corresponds to a section $s\in H^{0}(X,L)$ (defined
up to a non zero scalar). The Poincar\'e--Lelong equation is
$$dd^c\log\Vert s\Vert^2=\delta_D-c_1(L)$$
where $\delta_D$ is the Dirac distribution "integration on $D$" and
$c_1(L)$ is the first Chern form associated to the hermitian line
bundle $L$.

A divisor on $X$ is said to be {\it simple normal crossing} (snc for
short) if $D=\sum D_i$ with $D_i$ smooth and locally for the
Euclidean topology, we can find coordinates $x_1,\dots, x_n$ on $X$
for which $D_i=\{x_i=0\}$ and $D=\{x_1\cdots x_r=0\}$. If $D$ is
snc, we can introduce the sheaf of differentials on $X$ with
logarithmic poles along  $D$: $\Omega^1_X(\log(D))$; this is the
sheaf of meromorphic differentials $\omega$ which may locally be
written as $\omega=\sum_{i=1}^ra_i{{dx_i}\over{x_i}}+\alpha$ where
$a_i$ are $C^{\infty}$ functions and $\alpha$ is a smooth
differential. If $D$ is a divisor on a variety, we let $D_{red}$ be
the {\it reduced} divisor having the same support as $D$.

Suppose that $Y$ is a {\it compact} Riemann surface and $f:Y\to X$
is an analytic map such that $f(Y)\not\subset D$. Since $Y$ is
without boundary, Stokes theorem gives
$$\deg(f^\*(D))=\int_Yf^{\ast}(c_1(L)).$$
Thus {\it the degree} of the restriction of the divisor $D$ to $Y$
can be interpreted as the {\it area} of $Y$ with respect to the
measure defined with $c_1(L)$. We will denote by $N^{(1)}_D(Y)$ the
degree of $f^{\ast}(D)_{red}$; observe that $N^{(1)}_D(Y)=\sum_{z\in
Y}\min\{ 1,v_z(f^\*(D))\}$ thus, in the number fields--function
fields analogy, it corresponds to the radical defined in the
previous section. In the sequel we will denote by $(L,Y)$ the
integral number $\deg(f^\*(L))$ (omitting the reference to $f$ if
this is clear from the contest).

Suppose that $X=X_1\times X_2$ where $X_i$  are compact Riemann
surfaces and $p_i:X\to X_i$ are the projections. Suppose that $D_2$
is a reduced divisor on $X_2$ and $D:=p_2^\*(D_2)$. We consider
$(X;D)$ as an isotrivial family of curves with divisors over $X_2$
via $p_2$.

Let $Y$ be a compact Riemann surface with an analytic map $f:Y\to
X$. Call $f_i:=p_i\circ f$. For $i=1,2$, define $R_{f_i}$ as the
Ramification term: $R_{f_i}:=\sum_{z\in Y}(Ram(f_i)-1)$. The Hurwitz
formula gives
$f_2^\*(\Omega^1_{X_2}(\log(D_2)))\hookrightarrow\Omega^1_Y(f^\*(D)_{red})$.
Thus a double application of the Hurwitz formula gives
$$\deg(f^\*(K_{X/X_1}(D)))\leq N^{(1)}_D(Y)+R_{f_1}+\chi(X_1)\deg(f_1);$$
Which is the analogue (which holds with $\epsilon=0$) of the $abc$
conjecture over function fields in the isotrivial case.

When $f:Y\to X$ is an algebraic map between smooth projective
curves, we will denote by $[Y:X]$ the degree of the pull back, via
$f$, of a generic point; it coincides with the degree of the field
extension $\CC(Y)/\CC(X)$.

Suppose that $Y$ is not compact. In this case we suppose that $Y$ is
{\it parabolic} equipped with an exhaustion function $g$: an
exhaustion function is a unbounded function $g$ such that
$dd^c(g)=\delta_S$ where $S=\sum P_i$ is a reduced divisor of finite
degree and near each $P_i$ we can find an harmonic function $h_P$
such that $g=\log\vert z-P\vert^2+h$. Remark that $g$ is harmonic
outside $S$,

\noindent {\it Examples}: a) $\Bbb C$ with the function $\log\vert
z\vert^2$

\noindent b) If $\pi:Y\to \Bbb C$ is a proper map of degree
$[Y:\CC]$ (not ramified over $0$), then $(Y;\pi^\*(\log\vert
z\vert^2))$ is parabolic: thus every affine Riemann surface is
parabolic.

\noindent c) If $Y$ is parabolic and $E$ is a polar set in $Y$ then
$Y\setminus E$ is parabolic.

For more details on parabolic Riemann surfaces cf. [AS].

Parabolic Riemann surfaces are the ones where one can develop a
value distribution theory. We fix a parabolic Riemann surface
$(Y,g)$.

Suppose that $f:Y\to X$ is an analytic map. We define the {\it
intersection} of the hermitian line bundle $L$ with $Y$ as a {\it
function} on $\RR$:
$$(L;Y)(r):=\int_{-\infty}^{\log(r)}ds\int_{g\leq s}f^\*(c_1(L))=\int_0^{r}{{dt}\over{t}}\int_{g\leq\log(t)}f^\*(c_1(L))$$
(in value distribution theory this is denote as $T_f(L)$; we choose
this notation to stress the analogy with intersection theory). The
intersection $(L,Y)(r)$ can be seen as {\it the average of the areas
of the disks $g\leq s$} for $s\leq\log(r)$. Up to a constant
$(L;Y)(r)$ do not depend on the  choice of the metric on $L$.

We will define the {\it non integrated counting function} as
$n_D(s):=\sum_{g(z)<s} v_z(f^\*(D))$ and the {\it the non integrated
radical function} as $n^{(1)}_D(s):=\sum_{g(z)<s} \inf\{ 1,
v_z(f^\*(D))\}$; the n.i. counting function measure the growth of
the degree of the divisor $f^\*(D)$ on the disk $g<s$ and the n.i.
radical plays the role of the radical . We define the {\it
integrated counting function} and the {\it integrated radical} as
$$N_D(Y)(r):=\int_{-\infty}^{\log(r)}n(f^{\*}(D),s)ds \; {\rm and}
\;
N^{(1)}_D(Y)(r):=\int_{-\infty}^{\log(r)}n^{(1)}(f^{\*}(D),s)ds\;\;
{\rm resp.}$$

In the  same way we define a {\it non integrated characteristic
function or ramification term}: the form $\partial g$ is holomorphic
outside $S$ thus we define $r_g(s):=\sum_{g(z)<s}v_z(\partial g)$,
where the sum is extended to points not in $S$. For instance, if $Y$
is a proper covering of $\CC$, then $r_g(s)$ is the degree of the
part of the ramification divisor of the covering supported in $g<s$.
Thus we define the {\it integrated characteristic function} as
$$\chi(Y)(r):=\int_{-\infty}^{\log(r)}r_g(s)ds.$$

Integrating the Poincar\'e--Lelong equation we obtain the {\it First
Main Theorem} of Value distribution theory (cf. [NE] or [HA]):
Suppose that $X$ is smooth projective , $D$ and $L$ are as before
and $f:Y\to X$ then we can find an explicit constant $C$,
independent on $r$, such that
$$(L;Y)(r)=N_D(Y)(r)-\int_{g=r}\log\Vert s\Vert^2d^cg + C;$$
The term $m_D(Y,r):=-\int_{g=r}\log\Vert s\Vert^2d^cg$ is called
{\it the proximity function} and measure the average of the inverse
of the distance of the image of boundary of $g\leq r$ from $D$.

Suppose that $X:=X_1\times X_2$ and $D:=p_2^{\*}(D_2)$ as before.
Suppose that $f:Y\to X$ is an analytic map from a parabolic Riemann
surface. The {\it Second Main Theorem} of Value Distribution Theory
([NE] and [HA]) can be stated in this way:
$$(K_{X/X_1}(D);Y)(r)\leq
N_D^{(1)}(Y)(r)+\chi(Y)(r)+O(\log(r(K_{X/X_1}(D);Y)(r)))$$ where the
involved constant is independent on $r$. Thus,  the analogue of the
$abc$ conjecture in the isotrivial case is the second main theorem.

Let $Y$ be  parabolic, $B$ be a compact Riemann surface and $f:Y\to
B$ be an analytic map. Let $P\in B$ be a point and equip $\cO_B(P)$
with a smooth metric. In analogy with the algebraic case, we will
denote $[Y:B](r)$ the function $(\cO_B(P);Y)(r)$.

\endsection

\

\section Statement of the main theorems\par

\

In this section we will state the main theorems, namely the
$abc$--conjecture over function fields and make the first easy
reductions.

The object of study is a set $(X,D, B,p )$ where, $X$ is a smooth
projective surface, $D$ is a simple normal crossing divisor on $X$,
$B$ is a smooth projective curve and $p:X\to B$ is a non constant
morphism. We will also fix an ample line bundle $H$ equipped with a
smooth positive metric.

We will explain the proof of the two theorems below, one is in the
algebraic  and the other in the analytic setting. They correspond
each other in the analogy and we will see that the proofs are, {\it
mutatis mutandis}, very similar.

\label maitheoremalg. theorem\par\thm {\rm ($abc$ algebraic
version)} Let $p:X\to B$ and $D$ as above and $\epsilon>0$. Then,
given a smooth projective curve $Y$ and a morphism $f:Y\to X$ whose
image is not contained in $D$,  the following inequality holds
$$(K_X(D);Y)\leq
(1+\epsilon)(N_D^{(1)}(Y)+\chi(Y))+O_\epsilon([Y:B]).$$ The involved
implicit constants depend only on $X, D, p$ and $\epsilon$.\endthm

To avoid trivialities we supposed that the morphism $p\circ f:Y\to
B$ is non constant.

\label maintheoreman. theorem\par\thm {\rm ($abc$ analytic version)}
Let $p:X\to B$ and $D$ as above and $\epsilon>0$. Let $(Y,g)$ be a
parabolic Riemann surface and $f:Y\to X$ be a holomorphic map with
dense image. Then the following inequality holds
$$(K_X(D);Y)(r)\leq(1+\epsilon)(N_D^{(1)}(Y)(r)+\chi(Y)(r))+O_\epsilon([Y:B](r)+\log(r(H;Y)(r))). \;\; //$$
The involved implicit constants depend only on $X, D, p, f$ and
$\epsilon$ but independent on $r$.
\endthm

The symbol $//$ means that the inequality holds outside a set of
finite Lebesgue measure.

\ssection Reductions and observations\par a) The theorems remain
true if we change $X$ by a blow up (but the involved constants may
vary). Consequently they are statements about the algebraic curve
$X_K$ where $K$ is the function field $\CC(B)$. Moreover we may
suppose that every irreducible component of $D$ dominates $B$.

b) In order to prove the theorems we may take finite extensions of
the base field $K:=\CC(B)$: if the theorem is true over a finite
extension, it is true over it, and conversely. We may, and we will,
suppose for instance that $B$ is hyperbolic.

c) By the semistable reduction theorem, we may suppose that $K_X(D)$
is nef and big and that the fibres of $p$ are reduced simple normal
crossing. Incidentally, this shows that the hypothesis that $D$ is
simple normal crossing is unnecessary.

d) Suppose that $X=\PP_1\times B$ and $D$ is the pull back, via the
the first projection of the divisor $0+1+\infty$; then the theorems
give the "classical" {\it strong} $abc$--conjecture over function
fields.

e) Suppose that, in the analytic situation, $(Y,g)=(\CC, \log\vert
z\vert^2)$, $X=\PP_1\times \PP_1$, $D$ is the pull back via the
first projection of a divisor $\sum_{i=1}^d a_i$ and
$f:=(f_1,id):\CC\to X$, where $f_1$ is a meromorphic function, then
the theorem becomes (in the standard notation of Nevanlinna theory)
$$(d-2)T_{f_1}(r)\leq
(1+\epsilon)\sum_iN^{(1)}(a_i,f)+O_\epsilon(\log(rT_{f_1}(r)))\;\;\;
//.$$
which is essentially (up to the factor $\epsilon$) the Nevanlinna
Second Main Theorem.

\endssection

\endsection

\section The tautological inequality\par

\

Let $X$ be a smooth projective variety and $D\subset X$ be a simple
normal crossing divisor. Let $\Omega_{X/k}^1(\log(D))$ be the sheaf
of differentials of $X$ with logarithmic poles along $D$ and
$\pi:\PP:=\Proj(\Omega_{X/k}^1(\log(D)))\to X$. We will denote by
$\LL$ the tautological line bundle on $\PP$. We also fix an ample
line bundle $H$ on $X$

Let $Y$ be a smooth projective curve and $f\colon Y\to X$ a map. The
induced map $f^\*: f^\*(\Omega_{X/k}^1(\log(D)))\to
\Omega_{Y/k}^1(f^{\*}(D)_{red})$ induce a morphism $f'\colon Y\to
\PP$. By definition  we have an inclusion
$f'^\*(\LL)\hookrightarrow\Omega_{Y/k}^1(f^{\*}(D)_{red})$
consequently we obtain the {\it tautological inequality} \labelf
geotaut\par$$(\LL; Y)\leq \chi(Y)+N_D^1(Y)\eqno{{(\numfo)}
}$$\advance\ssnu by1 where $\chi(Y):=2g(Y)-2$ is the Euler
characteristic of $Y$.

Suppose now that $k=\CC$. We Suppose that $H$ is equipped with a
positive (Kh\"aler) metric $\omega$. We also equip $\cO(D)$  and
$\Omega_{X/k}^1(\log(D))$ with a  smooth metric. Remark that the
metric on $\Omega_{X/k}^1(\log(D))$ induces a metric on $\LL$.

Suppose that $(Y; g)$ is a parabolic Riemann surface and $f\colon
Y\to X$ is an analytic map whose image is not contained in $D$. The
tautological inequality is the analogue of the \ref{geotaut} in this
contest.

\label tautin. theorem\par\thm (Analytic tautological inequality)
Let $f\colon Y\to X$ as above and $f'\colon Y\to\PP$ the induced
map. The following inequality holds \labelf taut\par$$(\LL;Y)(r)\leq
N^1_D(Y)(r)+\chi(Y)(r)+ O(\log(H;Y)(r)) \;\;\; //\eqno{{(\numfo)}}$$
\endthm

The tautological inequality above is an important push forward in
the analogy between the algebraic and the analytic theory of maps of
Riemann surfaces in projective varieties. It is very important
because it translate the problems of defect type in Nevanlinna
theory to problems of geometrical nature: If one prove that some
intersection is upper bounded by the intersection with $\LL$, one
will deduce an inequality in the spirit of the Second Main Theorem
of Nevanlinna theory.

\rmk One may wonder how much of the proofs of $abc$ performed in the
function field case can be done in the arithmetic situation.
Unfortunately, in the arithmetic situation, the geometric
interpretation of the radical via the tautological inequality  is
missing: one do not know what is the arithmetic meaning of the
radical.

\endrmk

\Proof Write the divisor $D$ as $\sum_iD_i$. Locally on $X$ we can
find coordinates $X_1,\dots, X_n$ in such a way that, there is an
$r\in\{ 0,\dots n\}$ such that,  each $D_i$ is given by $\{X_i=0\}$
and $D=\{X_1\cdot\dots\cdot X_r=0\}$. Since $X$ is compact, can
choose a positive constant such that, the singular $(1,1)$ form
$$\omega^{sm}:=A\omega+\sum_i{{d\Vert D_i\Vert\wedge d^c\Vert
D_i\Vert}\over{\Vert D_i\Vert^2}}$$ induces a smooth hermitian
metric on $T_{X/k}(-\log(D))$. We introduce the singular $(1,1)$
form
$$\tilde\omega:=\omega+\sum_i{{d\Vert D_i\Vert\wedge d^c\Vert
D_i\Vert}\over{\Vert D_i\Vert^2(\log(\Vert D_i\Vert))^2}}.$$ The
form $\tilde\omega$ induces a singular hermitian form on
$T_{X/k}(-\log(D))$: if we write an element of $T_{X/k}(-\log(D))$
as $t=\sum_{i=1}^r a_i{{X_i\partial}\over{\partial
X_i}}+\sum_{j=r+1}^nb_j{{\partial}\over{\partial X_i}}$ then
$\tilde\omega(t,t)$ is comparable to $\sum_i\vert
a_iX_i\vert^2+\sum_j\vert b_j\vert^2+\sum_i{{\vert
a_i\vert^2}\over{\log\vert X_i\vert^2}}$.

Let $\PP'$ be the projective bundle
$\Proj(\cO_X\oplus\Omega_{X/k}^1(\log(D)))$ over  $X$, and let $\MM$
be the tautological line bundle over it. The surjection
$\cO_X\oplus\Omega_{X/k}^1(\log(D))\to\Omega_{X/k}^1(\log(D))$
induces an inclusion $\PP\hookrightarrow\PP'$ and
$\cO_{\PP'}(\PP)=\MM$. Observe that, locally on $X$, $\PP'$ is
isomorphic to $X\times\PP_n$ and, we may choose  homogeneous
coordinates $[z_0:\cdots: z_n]$  on $\PP_n$ for wich the divisor
$\PP$ is given by $z_0=0$.

On the other side the inclusion
$\Omega_{X/k}^1(\log(D))\to\cO_X\oplus\Omega_{X/k}^1(\log(D))$
induces a rational map $h:\PP'\dashrightarrow\PP$. Let $q:\tilde
Z\to\PP'$ be the blow up along the section given by the projection
$\cO_X\oplus\Omega_{X/k}^1(\log(D))\to\cO_X$; it resolves the
indeterminacy of $h$ and we obtain a morphism $p:\tilde Z\to\PP$. By
construction we obtain $p^\*(\LL)=q^\*(\MM)(-E)$ where $E$ is the
exceptional divisor of $\tilde Z$.

The form $\omega^{sm}$ induces positive metrics
$\Vert\cdot\Vert^{sm}_\MM$ on $\MM$ and $\Vert\cdot\Vert^{sm}_\LL$
on $\LL$. We denote by $c_1(\MM)^{sm}$ and $c_1(\LL)^{sm}$ the
corresponding singular first Chern forms. We put on $\cO_{\tilde
Z}(E)$ the metric for  which the isomorphism $p^\*(\LL)\simeq
q^\*(\MM)(-E)$ become an isometry. The form $\tilde\omega$ induces
singular metrics $\Vert\cdot\Vert_\MM^s$ on $\MM$ and
$\Vert\cdot\Vert_\LL^s$ on $\LL$. We denote by $c_1(\MM)^s$ and
$c_1(\LL)^s$ the corresponding singular first Chern forms.

The morphism
$$\eqalign{\cO_Y\oplus
f^\*(\Omega_{X}^1(\log(D)))&\longrightarrow\Omega^1_{Y}(S+f^{-1}(D))\cr
(a,\alpha)&\longrightarrow a\partial(g)+f^\*(\alpha)\cr}$$ induces
maps $f_1\colon Y\to\PP'$ and $\tilde f:Y\to\tilde Z$. By
construction $p\circ\tilde f=f'$.

Locally on $Y$, we can write $f$ as $(g_1,\cdots, g_n)$ and
$f_1=(g_1,\cdots, g_n)\times
[g':{{g_1'}\over{g_1}}:\cdots:{{g_r'}\over{g_r}}:g'_{r+1}:\cdots:g'_n]$.

We apply the first main theorem to the map $f_1$ the line bundle
$\MM$ equipped with $\Vert\cdot\Vert_\MM^s$ and the divisor $\PP$.
Remark that, even if the metric $\Vert\cdot\Vert_\MM^s$ is singular,
we can apply the FMT because it is locally integrable.

By the local computation of $f_1$ we see that
$N_\PP(Y)(r)=\chi(Y)(r)+N^{(1)}_D(Y,r)$. Thus we obtain
$$\int_{-\infty}^{\log(r)}dt\int_{g\leq
t}f_1^\*c_1(\MM)^s=\chi(Y)(r)+N^{(1)}_D(Y,r)-2\int_{g=r}\log\Vert\PP\Vert_\MM^s
d^cg.$$ Consequently we obtain
$$\eqalign{\int_{-\infty}^{\log(r)}dt\int_{g\leq
t}f'^\*(c_1(\LL)^s)&=
\cr=\chi(Y)(r)+N^{(1)}_D(Y,r)-&2\int_{g=r}\log\Vert\PP\Vert_\MM^s
d^cg-N_E(Y)(r)+2\int_{g=r}\log\Vert E\Vert d^cg+O(1).\cr}$$

We claim that $(\LL,Y)(r)$ is smaller then
$\int_{-\infty}^{\log(r)}dt\int_{g\leq
t}f'^\*(c_1(\LL)^s)+O(\log(H;Y)(r))$. The  intersection $(\LL,Y)(r)$
can be computed using the metric $\Vert\cdot\Vert^{sm}_\LL$. Outside
$D$, we can find a function $h$ such that
$\Vert\cdot\Vert^{sm}_\LL\cdot h=\Vert\cdot\Vert^{s}_\LL$. Thus
$c_1(\LL)^s=c_1(\LL)^{sm}-dd^c\log(h)$. Computing the two metrics
locally, again by compactness of $X$, we obtain that $h\ll
\prod\log^2\Vert D_i\Vert$.   Consequently
$$\int_{-\infty}^{\log(r)}dt\int_{g\leq
t}f'^\*(c_1(\LL)^s)=(\LL;Y)(r)-\int_{\infty}^{\log(r)}dt\int_{g\leq
t}dd^c\log(h)+O(1);$$ by applying Stokes Theorem twice, we find that
$$\int_{\infty}^{\log(r)}dt\int_{g\leq
t}dd^c\log(h)=\int_{g=\log r}\log(h)d^cg+O(1).$$ This last term can
be bounded as follows
$$\eqalign{\int_{g=\log r}\log(h)d^cg\ll&\sum_i\int_{g=\log r}\log(\log^2\Vert
D_i\Vert)d^c g\cr &\leq\sum_i2\log\int_{g=\log r}\vert\log\Vert
D_i\Vert\vert d^c g\cr
&\leq\sum_i2\log((\cO_X(D_i);Y)(r))\ll\log((H;Y)(r)).\cr}$$ The
claim follows.

We compute now, locally, $\tilde f^\*(\Vert E\Vert)(z)$ and
$f_1^\*(\Vert\PP\Vert_\MM^s)(z)$. Let $z$ be a local coordinate on
$Y$ and let $\partial_z$ be the corresponding local generator of the
tangent bundle of $Y$. Define
$\tilde\omega(z):=f^\*(\tilde\omega)(\partial_z)$. The local
expression of $f$, $f_1$ etc. implies that $\tilde f^\*(\Vert
E\Vert)^2(z)={{\tilde\omega(z)}\over{\tilde\omega(z)+\vert\partial_z
g\vert^2}}$ and $\tilde
f^\*(\Vert\PP\Vert_\MM^s)^2(z)={{\vert\partial
g\vert^2}\over{\tilde\omega(z)+\vert\partial_z g\vert^2}}$.

In order to conclude, we need to find an upper bound for
$$T(r):=\int_{g=r}\log{{\tilde\omega(z)}\over{\vert\partial_z g\vert^2}}d^cg.$$

We can find a function $F$ such that $f^\*(\tilde\omega)=Fdg\wedge
d^cg$. The function $F$ will be $\vert z\vert^2\times {\rm smooth}$
in the neighborhood of the poles of $g$ and in general
$F(z)={{\tilde\omega(z)}\over{\vert\partial_z g\vert^2}}$.

Let $$S(r):=\int_{\infty}^{\log(r)}dt\int_{g\leq
t}f^\ast(\tilde\omega).$$ Fubini Theorem gives
$$S(r)'=\int_{-\infty}^{\log(r)}dt\int_{g=t}Fd^cg;$$ Thus, cancavity
of the $\log$ gives $$\log(S^{(2)}(r))\geq T(r).$$

The following lemma is well known and elementary (for a proof cf.
[GK])

\label lang2. lemma\par\lemma Let $H$ be a derivable positive
increasing function. For every positive $\epsilon$, there exists a
subset $E\subset\Bbb R$ with
$meas(E)\leq\int_{1+\epsilon}^{\infty}{{1}\over{x\log^{1+\epsilon}(x)}}dx<\infty$,
such that, for every $x\not\in E$,
$$H'(x)\leq H(x)\log^{1+\epsilon}(H(x)).$$
\endlemma

We apply Lemma \ref{lang2} twice and we find that
$$T(r)\leq\log(
S(r)\log^{1+\epsilon}(S(r))\log^{1+\epsilon}(S(r)\log^{1+\epsilon}(S(r)))).$$

We will conclude if we find an upper bound for $S(r)$.

The following equality holds
$$-dd^c\log(\log^2(\Vert D_i\Vert))={{d\Vert D_i\Vert\wedge d^c\Vert
D_i\Vert}\over{\Vert D_i\Vert^2(\log\Vert
D_i\Vert)^2}}+{{1}\over{\vert\log(\Vert D_i\Vert)\vert}}\cdot
c_1(\cO(D_i)); $$ by the compactness of $X$, the last term on the
right hand side is uniformly bounded; thus we can find a constant
$A$ such that
$$\tilde\omega\leq A\omega -\sum_idd^c\log(\log^2(\Vert D_i\Vert)).$$
Thus, again by applying Stokes,
$$S(r)\ll (H;Y)(r)-2\int_{g=r}\log(\log^2(\Vert D_i\Vert))d^cg
+O(1).$$ Since we can suppose that $\Vert D_i\Vert<\epsilon$ we
conclude.
\endsection

\

\section Currents associated to families of curves\par

\

Let $X$ be a projective variety (reduced and irreducible) and $H$ an
ample line bundle equipped with a smooth positive metric. We will
now show how to associate a closed positive current to the situation
we are interested in.

In the analytic situation we start with a map from a parabolic
Riemann surface to $X$ and the diophantine statement we are
interested in upper bounds {\it uniform in $r$}. Roughly speaking we
have maps from the Riemann surfaces $\{ z\in Y \; / \;
g(z)<\log(r)\}$ and we look for uniform upper bounds for their areas
(or better: the average of the areas over them) with respect to some
hermitian line bundle, in terms of their Euler characteristic and
the (set theoretical) intersection with the divisor at infinity.
Thus we can take a sequence of $r$'s which goes to infinity and do
not satisfy the wanted inequality and eventually find a
contradiction.

Similarly, in the algebraic case, we want to give uniform upper
bounds of the intersection (the height!) of the closed curves in $X$
in terms of the Euler characteristic of their normalization and the
(set theoretic) intersection with the divisor at infinity. Again we
take a sequence of smooth projective curves which do not satisfy the
inequality and find a contradiction.

We will show now, that in both situations we can associate to the
involved sequences a closed positive current $T$ on $X$ (and on
other varieties constructed during the proof). The proof of the
theorem will work with the properties of $T$ and only at the end,
the definition of it will give the statement in the analytic or in
the geometric case.

We will consider two situations:

\smallskip

a) {\it The analytic situation $S^{an}$}: A parabolic Riemann
surface and a holomorphic map $f\colon Y\to X$.

\smallskip

b){\it The algebraic situation $S^{alg}$}: A sequence of smooth
projective curves $\{Y_n\}_{n\in\Bbb N}$ and algebraic maps $f_n:
Y_n\to X$.

\smallskip

Before we start the construction, we have to show that the currents
associated to $S^{an}$, even if they are not closed, they are
"closed enough".

Let $(Y,g)$ be  a parabolic Riemann surface equipped with a positive
singularity. Let $f:Y\to X$ be a holomorphic map.

We show now that the intersection of $Y$ with exact forms is
essentially irrelevant.

\label closednull. lemma\par\lemma Let $1/2>\epsilon>0$. Let
$\alpha$ be a smooth exact $(1,1)$ form on $X$ then
$$(\alpha, Y)(r)=O_\epsilon(((H;Y)(r))^{{1/2}}\log^{1-\epsilon}(r(H;Y)(r))) \;\; //$$
where the involved constants depend only on $\epsilon$.\endlemma

\Proof Since $(\cdot, Y)(r)$ is a positive current, it will suffice
to prove the theorem when $\alpha$ is $\overline\partial\beta$ for a
smooth $(1,0)$ form $\beta$.

By Stokes theorem
$$\eqalign{(\alpha; Y)(r)&=\int_{0}^{r}{{dt}\over{t}}\int_{g=t}\beta\cr
&= \int_{-\infty<g\leq\log(r)}dg\wedge\beta.\cr}$$

Since $X$ is compact, we have that $\beta\wedge\overline\beta\ll
c_1(H)$. Consequently, Cauchy--Schwartz inequality gives
$$\left\vert\int_{-\infty<g\leq\log(r)}dg\wedge\beta\right\vert
\leq
2\pi\left\vert\int_{-\infty<g\leq\log(r)}\beta\wedge\overline\beta\right\vert^{1/2}\cdot
\left\vert\int_{-\infty<g\leq\log(r)}dg\wedge
d^cg\right\vert^{1/2}.$$ We apply again \ref{lang2} and we obtain
that, outside a set of finite Lebesgue misure,
$$\left\vert\int_{-\infty<g\leq\log(r)}\beta\wedge\overline\beta\right\vert\ll
(H,Y)(r)\log^{1+\epsilon}((H;Y)(r)).$$ Since (again by Stokes)
$\int_{-\infty<g\leq\log(r)}dg\wedge d^cg\leq C\log(r)$, for a
suitable $C$, we conclude.

\label algebraicremark1. remark\par\rmk  In the algebraic setting,
if $Y$ is a smooth algebraic curve and $f\colon Y\to X$ is a map,
the current $\alpha\to\int_Yf^{\ast}(\alpha)$ is closed by Stokes
theorem. We observe that the lemma above tells us that, up a
negligible term, the current $(\cdot;Y)(r)$ is closed. Consequently,
up to this negligible term, this is another analogy between the two
situations, we will see now how to push forward this.
\endrmk

Let $(X,H)$ as above. In the analytic situation consider the set of
currents
$$\eqalign{ T_r&: A^{1,1}(X)\longrightarrow \RR\cr
&\alpha\longrightarrow
{{1}\over{(H;Y)(r)}}\int_{0}^r{{dt}\over{t}}\int_{g<t}f^\*(\alpha);\cr}$$
in the geometric situation consider the set of currents
$$\eqalign{T_n&: A^{1,1}(X)\longrightarrow\RR\cr &\alpha\longrightarrow
{{\int_{Y_n}f_n^\*(\alpha)}\over{(H;Y_n)}}.\cr}$$

\smallskip

In both situations they are families of positive currents bounded
for the standard norm on $A^{1,1}$ (but also for the $L^{\infty}$
norm); consequently we can extract from them a sequence converging,
in the weak topology, to a positive current $T$.

In the algebraic situation the current $T$ is closed because it is
limit of closed currents.

In the analytic case, due to Lemma \ref{closednull}, we can choose
the sequence in order to obtain a sequence $T_{r_n}$ such that
$dT_{r_n}\to 0$. Observe that, even if the involved map $f$ is
algebraic, the height of $r$ at least  $A\log(r)$ (for a suitable A)
consequently the Lemma apply.

\label current. definition\par\defi The closed positive current $T$
constructed above,  will be called {\rm the current} associated to
the (geometric or analytic) situation.
\enddefi

Observe that $T$ is non zero because $T(c_1(H))=1$.

\rmk the article {\it the} in the definition is not completely
correct. Indeed the current $T$ depends on the choice of the
subsequences involved. The reader will check that we will use only
properties which hold for {\it every} sequence as above.\endrmk

\label interdivisor. remark\par\rmk (Important) Since the current
$T$ is closed, we can unambiguously compute it on a Cartier divisor
$R$ of $X$: put an arbitrary smooth metric on $\cO_X(R)$ and define
$T(R):=T(c_1(\cO_X(R))$; this number do not depend on the chosen
metric. Of course, since $T$ is positive, if $R$ is ample (resp.
nef) then $T(R)>0$ (resp. $T(R)\geq 0$). One should see $T$ as a
class in the dual of the positive cone of $NS(X)_\RR$ and interpret
$T(R)$ as an intersection number.
\endrmk

\endsection

\

\section First approach to the theorems\par

\

In this section we will explain the proof by McQuillan of theorems
\ref{maitheoremalg} and \ref{maintheoreman}. We recall that we are
assuming the reductions made after the statements. We begin by
fixing some notations: $\pi:\PP:=\Proj(\Omega_{X}^1(\log(D)))\to X$
will be the projective bundle associated to the sheaf of
differentials with logarithmic poles around $D$. Let $\LL$ be the
tautological bundle over $\PP$. We denote by $F$ a smooth fibre of
$p$. We fix an ample divisor $H$ on $X$; observe that we may suppose
that there exists $\epsilon'$ such that, for every curve $Y$ not
contained in a fibre, $(Y,K_X(D))\geq\epsilon' (Y;H)$ (and similarly
in the analytic setting).

The theorems are proved  by contradiction.

a) {\it Algebraic situation}: We suppose that there is a sequence of
smooth projective curves and morphisms $f_n:Y_n\to X$ such that
\labelf
abshypothesisalg\par$$\lim_{n\to\infty}{{N_D^{(1)}(Y_n)+\chi(Y_n)}\over{(H;Y_n)}}<\infty
\;\;\;\; {\rm and }\; \;\;\;
\lim_{n\to\infty}{{[Y_n:B]}\over{(H;Y_n)}}=0.\eqno{{(\numfo)}}$$\advance\ssnu
by1 Consequently we can construct a closed and positive current $T$
on $X$ associated to the sequence. We can rise each map $f_n$ to a
map $f'_n:Y_n\to \PP$. Each $f_n'$ give rise to a closed positive
current
$$\eqalign{T'_n:&A^{(1,1)}(\PP)\longrightarrow \RR\cr
\alpha&\longrightarrow
{{1}\over{(H;Y_n)}}\int_{Y_n}f_n'^\ast(\alpha).\cr}$$

{\it Because of the hypothesis \ref{abshypothesisalg}}, we can
extract from the sequence above a subsequence converging to a closed
positive current $T'$ on $\PP$. By construction we have that
$\pi_\ast(T')=T$.

The theorem will be proved if we show that
$$T'(\LL-\pi^\ast(K_{X/D}(D)))\geq 0.$$

\smallskip

b) {\it Analytic situation}: We suppose that there exists a sequence
of real numbers $r_n$ such that \labelf
abshypothesisan\par$$\lim_{n\to\infty}{{N_D^{(1)}(Y)(r_n)+\chi(Y)(r_n)}\over{(H,Y)(r_n)}}<\infty\;\;\;\;{\rm
and }\;\;\;\;\lim_{n\to
\infty}{{[Y;B](r_n)}\over{(H;Y)(r_n)}}=0.\eqno{{(\numfo)}}$$\advance\ssnu
by1 Thus, again we can associate to this sequence the closed
positive current $T$ on $X$. We can rise the map $f$ to a map
$f':Y\to\PP$. Each $r_n$ give rise to a closed positive current
$$\eqalign{T'_n:&A^{(1,1)}(\PP)\longrightarrow \RR\cr
\alpha&\longrightarrow
{{1}\over{(H;Y)(r_n)}}\int_{0}^{r_n}{{dt}\over{t}}\int_{g\leq
t}f'^\ast(\alpha).\cr}$$ {\it Because of the hypothesis
\ref{abshypothesisan}}, we can extract from the sequence above a
subsequence converging to a closed positive current $T'$ on $\PP$.
By Lemma \ref{closednull}, if we choose suitably the subsequence,
the currents will converge to a closed current because we can
suppose that $\LL +\pi^\*(H)$ is ample.  By construction we have
that $\pi_\ast(T')=T$.

Again, the theorem will be proved if we show that
$$T'(\LL-\pi^\ast(K_{X/D}(D)))\geq 0.$$

\label vertical. remark\par\rmk Observe that  \ref{abshypothesisalg}
and \ref{abshypothesisan} imply that if $V$ is a vertical divisor on
$X$ (the map $p\vert_V:V\to B$ is not dominant) then $T(V)=0$.
\endrmk

Let  $p:X\to B$ the given fibration; $S:=\{ P_i\}$ be the set of
singular points of the fibers of $p$. By the reduction assumed, we
may suppose that  each $P_i$ belongs to exactly two smooth component
of a fibre which intersect properly on $P_i$; denote these two
components $B^i$ and $C^i$. We will denote by $I_S$ the ideal sheaf
of the subscheme defined by $S$ on $X$. The fibration $p$ give rise
to an exact sequence \labelf exactseq1\par$$0\to
p^\ast(\Omega^1_B)\longrightarrow\Omega^1_X(\log(D))\longrightarrow
I_S\otimes K_{X/B}(D)\to 0.\eqno{{(\numfo)}}$$\advance\ssnu by1 Let
$b:\tilde X\to X$ be the surface obtained taking the blow up of $X$
over the $P_i$'s and $E:=\sum_iE_i$ be the exceptional divisor
($E_i$ being the exceptional divisor over $P_i$). The exact sequence
\ref{exactseq1} give rise to an injection $\iota:\tilde X\to\PP$; by
construction
$$\iota^\*(\LL)=b^\*(K_{X/B}(D))(-E).$$
Denote by $\Delta$ the image of $\iota$; it is  a divisor on  $\PP$.
Observe that the fibre of $\pi$ over the $P_i$'s is contained in
$\Delta$. Denote by $U$ the open set $\PP\setminus\Delta$.

One of the main tools of the proof is the

\label pminusdelta. proposition\par\prop We have that
$$\pi_\ast(\Bbb I_UT')=0.$$
\endprop

The proposition above is not a surprise! Essentially it tells us the
following: The area of a curve on $X$ has a vertical and an
horizontal component (with respect to $p$); If the theorems are
false, then  we can find a sequence of curves  for with the vertical
component of the area is much bigger then the horizontal one. Thus
the limit of the curves must be vertical.

\Proof We fix  Kh\"aler forms $\omega$ on $X$ and $\eta$ on $B$. In
order to prove the proposition it suffices to prove the following:
Let $V$ be an open set of $\PP$ such that $\overline
V\cap\Delta=\emptyset$ ($\overline V$ being the closure of $V$ in
the Euclidean topology), then $\pi_\*(\Bbb I_VT)=0$. To prove this
we will show the following: there exists a constant $A_V$ (depending
on $V$ and the metrics) such that the following holds: if $W$ is an
open Riemann surface and $h:W\to X$ is an holomorphic map such that
the image of $h':W\to\PP$ is contained in $V$, then
$h^\*(\omega)\leq A_V(p\circ h)^\*(\eta)$. The conclusion follows
because $T$ (resp. $T'$) is zero on the fibres of $p$ (resp. on the
fibres of $p\circ\pi$).

Fix such a $V$, Observe that $\pi(\overline V)$ is a closed set of
$X$ which do not contains the singular points of the fibers $P_i$.
By compactness of $X$, we can cover $\pi(\overline V)$ by a finite
set of disks $B_j$ not containing the $P_i$'s. We may restrict our
attention to each of the $B_j$: thus we may suppose that:

\noindent -- $X=\{ (z,w)\in\CC^2\; /\; \vert z\vert<1 \; \vert
w\vert<1\}$, $B=\{ z\in C\; /\vert z\vert<1\}$ and $p(z,w)=z$;

\noindent -- $\omega= \sqrt{-1}(dz\wedge d\overline z+dw\wedge
d\overline w)$ and $\eta=\sqrt{-1}(dz\wedge d\overline z)$;

\noindent -- $D=\{ w=0\}$ and the exact sequence \ref{exactseq1} is
the split exact sequence
$$0\to \cO_X
dz\longrightarrow\cO_Xdz\oplus\cO_X{{dw}\over{w}}\longrightarrow\cO_X{{dw}\over{w}}\to
0;$$

\noindent -- consequently $\PP=X\times\PP^1$ and
$\Delta=X\times\{[0:1]\}$; we may then suppose that there exists a
positive constant $a$ such that $V\subseteq \{ (z,w)\times[x:y] \;
/\; \vert x\vert^2 >a\vert y\vert^2\}$.

\noindent -- $W:=\{ z\; /\; \vert z\vert<1\}$ and $h(z)=(h_1;h_2)$
and $h'(z)=(h_1;h_2)\times[h'_1:{{h'_2}\over{h_2}}]$.

The image of $W$ via  $h'$ is contained in $V$, we have that $\vert
h_1'(z)\vert^2>a\left\vert{{h'_2}\over{h_2}}\right\vert^2$. Thus
${{\vert h_2'(z)\vert^2}\over{\vert h_1'(z)\vert^2}}<{{1}\over{a}}$.
Since $h^{\*}(\omega)=\sqrt{-1}(\vert h_1'\vert^2+\vert
h'_2\vert^2)dz\wedge d\overline z$ and $h^{\*}(\eta)=\sqrt{-1}(\vert
h_1'\vert^2)dz\wedge d\overline z$ the proposition follows.

\smallskip

Since the line bundle $\LL$ is nef, as far as we are interested, we
may work as if the current $T$ is supported on $\Delta$. Indeed, by
the proposition,  $\Bbb I_UT=T-\Bbb I_\Delta T$ is a current which
intersect positively $\Bbb L$ and whose intersection with
$\pi^\*(K_X(D))$ is zero. Moreover, again, as far as we are
interested, via the proposition below, we can even suppose that it
is the push forward of a current on $\Delta$.

\label extensionofcurrents. proposition\par\prop Let $X$ be a smooth
variety and $E$ be a smooth divisor on it. Let $T$ be a closed
positive current of type $(1,1)$ over $X$. Let $\iota\colon E\to X$
be the inclusion. Then there is a current $S$ on $E$ such that
$$\II_E\cdot T=\iota_\*(S).$$
\endprop

It is evident that $S$ is positive and closed.

\Proof In order to prove the proposition, we need to show that, if
$\omega$ is a form such that $\iota^{\ast}(\omega)=0$ then
$\II_E\cdot T(\omega)=0.$

Locally on $X$ we can suppose that $z_1,\dots ,z_n$ are coordinates
on $X$ and $E=\{ z_n=0\}$. The theorem is proved if we show that for
every $i$ and for every $C^{\infty}$ function $a$ with compact
support, we have that $\II_E\cdot T(a\cdot dz_n\wedge d{\overline
z}_i)=\II_E\cdot T(a\cdot dz_i\wedge d{\overline z}_n)=0$.

Let $\psi$ be a positive function with compact support which is $1$
on the support of $a$. Since $\II_E\cdot T$ is positive, the
Cauchy--Schwartz inequality gives
$$\vert\II_E\cdot T(a\cdot dz_n\wedge d{\overline
z}_i)\vert^2\leq\vert\II_E\cdot T(\psi\cdot \sqrt{-1}dz_n\wedge
d{\overline z}_n)\vert^2\cdot\vert\II_E\cdot T(a\cdot
\sqrt{-1}dz_i\wedge d{\overline z}_i)\vert^2;$$ consequently  it
suffices to show that $\vert\II_E\cdot T(\psi\cdot
\sqrt{-1}dz_n\wedge d{\overline z}_n)\vert^2=0$. Since $T$ is of
type $(1,1)$ and closed, $\II_E\cdot T(\psi dd^c(\vert
z_n\vert^2)=\II_E\cdot T(\vert z\vert^2 dd^c\psi)=0$. But since
$dd^c\vert z_n\vert^2={{\sqrt{-1}}\over{2\pi}}dz_n\wedge d{\overline
z}_n$; the conclusion follows.

\smallskip

We apply the proposition above with $X=\PP$ and $E=\Delta$. Thus,
there is a closed positive current $S$ on $\Delta$ such that
$\II_\Delta T'=\iota_\*(S)$. Observe that, by functoriality,
$b_\*(S)=T$.

The proof of the theorem will be achieved if we can prove that
$S(-E_i)\geq 0$. In particular, if $X$ is smooth over $B$ and the
divisor $D$ is \'etale over $B$ then the theorem is proved; for
instance the isotrivial case (Nevanlinna Second Main Theorem) is
proved.

We reduced the difficulty to a local problem around the singular
points of the fibres. Most of the remaining of the proof will be of
local nature, but will will notice that one main point will be of
global nature.

\smallskip

The proof proceed by working on coverings of $X$; Let $Q_1,\dots,
Q_r$ be the points of $B$ where $p$ is not smooth. we fix another
point $Q$ on $B$. For every $m$ sufficiently big, we can find a
covering $B_m\to B$ which is cyclic of order $m$, totally ramified
over $Q_1,\dots, Q_r$ and \'etale over
$B\setminus\{Q,Q_1,\dots,Q_r\}$.

In the algebraic situation, the normalization $Y_{n,m}$ of the
curves $Y_n\times_BB_m$ are such that ${{\chi(Y_{n,m})}\over{m}}\leq
\chi(Y_n)+A$ (with $A$ independent on $m$).

In the analytic situation, the normalization $Y_m$ of $Y\times_BB_m$
is a parabolic Riemann surface, with as positive singularity, the
pull back of $g$ (which we will denote by $g_m$). Also in this
situation ${{\chi(Y_{m})(r)}\over{m}}\leq +\chi(Y)(r)+A$.

Consider the surface $X_m:=X\times_BB_m\buildrel{g_m}\over\to X$.
Let $p_m:X_m\to B_m$ be the second projection.

The surface $X_m$ is  normal and $\Bbb Q$--factorial. Denote by
$D_m$ the divisor $g_m^\ast(D)$. For every $i$ there is a unique
singular point $R_i$ over the $P_i$. For every $i$, there is an
analytic neighborhood of $R_i$ isomorphic to an analytic
neighborhood of the surface $Z^m=XY$ ($R_i$ corresponds to
$(0,0,0)$). Denote by $U_m$ the open surface $X_m\setminus\{R_i\}$.

\smallskip

In the next subsection we will justify the following properties:

\noindent -- There exists a normal $\Bbb Q$ factorial variety
$\PP_m$ with a $\Bbb Q$ line bundle $\LL_m$ and a projection
$\pi_m:\PP_m\to X_m$. Over $U_m$ there is a canonical isomorphism
$i_m:\PP_m\vert_{U_m}\leftrightarrow\PP(\Omega_{U_m}(\log(D_m))$
whose pull back of the tautological line bundle is the restriction
of $\LL_m$.

\noindent -- The analogue of the exact sequence \ref{exactseq1} over
$U_m$ give rise to an inclusion $U_m\to\PP(\Omega_{U_m}(\log(D_m))$;
Let $\tilde X_m$ be the closure of the image. Let $h_m:\tilde X_m\to
X_m$ be the projection and $\iota_m:\tilde X_m\to\PP_m$ be the
inclusion. The surface $\tilde X_m$ is again normal and $\Bbb
Q$--factorial. The preimage of $R_i$ is a Weil divisor and coincide
with the fibre over $R_i$ of $\pi_m$. Denote it by $E^m_i$. Moreover
$\iota_m^\*(\LL_m)=g_m^\*(K_{X_m/B_m}(D_m))(-\sum_iE_i^m)$ (this is
an equality between $\Bbb Q$--divisors). Denote by $V_m$ the open
set $\PP_m\setminus\iota(\tilde X_m)$ (remark that $V_m$ is smooth).

\noindent -- For every $R_i$ let $B^i_m$ and $C^i_m$ be the two
components of the fibre of $p_m$ meeting on $R_m$. We have that
$h_m^\*(B^i_m)=\tilde B^i_m+E_i^m$ and $h^\*_m(C^i_m)=\tilde
C^i_m+E^i_m$ where $\tilde B^i_m$ is birational to $B^i_m$ and
$\tilde C_i^m$ to $C_m^i$ respectively.

\noindent -- Let $f_{n,m}: Y_{n,m}\to X_m$ (resp. $f_m:Y_m\to X_m$
in the analytic situation); taking, if necessary, a subsequence of
the curves $Y_{n,m}$ (resp. of the $r_n$), we can construct, as
before, a closed positive current $T_m$ on $X_m$ (resp. $T'_m$ on
$\PP_m$) such that $g_{m,\*}(T_m)=T$; observe that we have to
normalize dividing by $m$. The value of the currents $T_m$ and
$T'_m$ on the fibres of $p_m$ and of $\pi_m\circ p_m$ respectively,
is  zero.

\noindent -- We can find a constant $A_m$ (depending on $m$) such
that, in the algebraic situation
$${{(\LL_n, Y_{n,m})}\over{m}}\leq N_D^{1}(Y_n)+\chi(Y_n)+A_m[Y_n;
B];$$ and in the analytic situation
$${{(\LL_n, Y_{m})(r)}\over{m}}\leq N_D^{1}(Y_n)(r)+\chi(Y_n)(r)+A_m[Y;
B](r);$$ thus the theorem will be proved if we show that there
exists a constant $A$ (independent on $m$) such that
$$T'_m(\LL_m-\pi_m^\*(K_{X_m/B_m}(D_m)))\geq {{A}\over{m}}.$$

Since the singular points of $\PP_m$ are contained in the image of
$\iota_m$; we can prove, exactly as before, that
$$\pi_{m,\*}(\II_{V_m}T')=0.$$ Since $\PP_m$ is $\Bbb Q$--factorial and
$\iota_m(\tilde X_m)$ is a divisor, the analogue of
\ref{extensionofcurrents} holds; thus there is a current $S_m$ on
$\tilde X_m$ such that $\II_{\iota_m(\tilde
X_m)}T_m'=\iota_\*(S_m)$. Moreover $h_{m,\*}(S_m)=T_m$.

The theorem is proved if we show that there is a constant $A$ such
that
$$S_m(-E^m_i)\geq {{A}\over{m}}.$$

Computing on the smooth part of $X_m$, we find that
$g_m^\*(B^i)=mB^i_m$ and $g_m^\*(C^i)=mC^i_m$.

On $\tilde X_m$, since $g_m^\*(B^i)=mB^i_m$ and $(g_m\circ
h_m)_\*(S_m)=T$ we have that  $S_m(h_m^\*(B_m^i))=0$ (cf. remark
\ref{vertical}). Thus $S_m(-E^i_m)=S_m(\tilde B_m^i)$.

\label btildecurrent. lemma\par\lemma Let $B$ be an effective
divisor on a projective  variety $X$ and $R$ be a closed positive
current on $X$ such that $\II_BR=0$, then $R(B)\geq 0$
\endlemma

We will postpone the proof of the lemma in the next subsection.

Because of the lemma, applied to $\II_{\tilde X_m\setminus \tilde
B_m^i}S_m$, we have that $S_m(\tilde B_m^i)\geq \II_{\tilde
B_m^i}S_m(\tilde B_m^i)$. Since $\tilde B_m^i$ and $\tilde C_m^i$
are disjoint, $\II_{\tilde B_m^i}S_m(\tilde B_m^i)=\II_{\tilde
B_m^i}S_m(\tilde B_m^i-\tilde C_m^i)$.

The divisor $\tilde B_m^i-\tilde C_m^i$ is $(h_m\circ g_m)^\*(
{{B^i- C^i}\over{m}})$, thus
$$S_m(-E^m_i)\geq {{(h_m\circ g_m)_\*(\II_{\tilde
B^i_m}S_m)(B_i-C_i)}\over{m}}.$$ Since $(h_m\circ g_m)_\*(S_m)=T$,
the following easy remark, applied to the couples $(T ;(h_m\circ
g_m)_\*(\II_{\tilde B_m^i}S_m)$ and $(T; (h_m\circ
g_m)_\*(\II_{\tilde X_m\setminus \tilde B_m^i}S_m))$, allows to
conclude:

Let $C$ be a divisor on a variety $X$ and $T$ a closed positive
current on $X$; then there exists a constant $A$ depending only on
$C$ and $T$ for which the following holds: for every closed positive
current $T_1$ with $T\geq T_1$, we have that $T(C)\geq T_1(C)+A$
(proof: take an ample bundle $H$ such that $C+H$ is ample and
compute $T$ and $T_1$ on $C+H$).

\ssection Extension of some results to singular varieties\par In
this subsection, we will explain how to extend the results we need
to the singular varieties involved in the proof. As explained in the
introduction, a natural approach to the proof is via integration on
stacks. Unfortunately, even for this easy orbifold case, we need to
develop the entire theory of integration on stacks; this is why we
prefer to deal with singular varieties.

A systematic approach to the tautological inequality and the other
needed results may be quite complicate, in particular it is not easy
to find the minimal hypotheses. This is why we develop just what we
need. Moreover this subsection will be sketchy.

{\it Metrized line bundles on singular varieties}: Let $X$ be a
reduced irreducible projective variety. Let $\cal L$ be a line
bundle on it. A continuous metric on $\cal L$ is a metric on the
fibres of it which varies continuously for the Euclidean topology.
We will say that a metric is {\it smooth} if locally we can embed
$X$ in a smooth variety $W$, $\cal L$ is the restriction of a line
bundle $\cal L_W$ on $W$ and the metric is the restriction of a
smooth metric on $\cal L_W$. We see that this is equivalent to ask
that, for every smooth variety $Y$ and map $f:Y\to X$, the induced
metric on the line bundle $f^\*(\cal L)$ is smooth. A (local)
section of $\cal L$ is said to be smooth if, locally it is the
restriction of a section on a smooth variety.

Observe that the sheaf $\Omega^{1,1}_X$ has a meaning on $X$:
$\Omega^1_X$ exists, and $\overline\Omega^1_X$ is its conjugate;
thus $\Omega^{1,1}:=\Omega^1_X\otimes\overline\Omega^1_X$. A $(1,1)$
form is said to be {\it smooth} if, locally it is the restriction of
a smooth form of a smooth variety. Similarly for functions.

Every line bundle on $X$ is difference of very ample line bundles,
thus every line bundle on $X$ admits a smooth metric.

Given a line bundle $\cal L$ on $X$ equipped with a smooth metric,
we can define its first Chern form in the following way: take a
(local) smooth section $f$ and $c_1(L):=-dd^c\log\Vert f\Vert^2$
outside the zeroes of $f$. Observe that $dd^c$ is well defined on
smooth functions and that $c_1(L)$ is a smooth $(1,1)$ form on $X$.
If we change the metric on $\cal L$ by another smooth metric, the
first Chern form varies by the $dd^c$ of a smooth function on $X$.

We gave examples to show that we can define all the objects we need
as restriction of similar objects defined over smooth varieties: in
particular we can define also the currents on $X$ and we can give a
meaning to closed and positive currents.

{\it Construction of $\PP_m$ and related objects}: The surface $X_m$
is smooth except on the points $R_i$. Near the $R_i$ it is
isomorphic to the surface $Z^m=XY$. Let $D_{\zeta, \xi}:=\{
(\zeta,\xi)\; / \vert \zeta\vert<1 ;\; \vert \xi\vert<1\}$. Let
 $\mu_n$ the cyclic group of the $m$--roots of the unity and let
$\theta_m$ be a generator of it; it acts on $D_{\zeta, \xi}$ with
the action $\zeta\to \theta_m\zeta$ and $\xi\to\theta_m^{-1}\xi$.
For every $i$, there is a neighborhood $V_i$ of the singular point
$R_i$ on $X_m$, isomorphic to $D_{\zeta, \xi}/\mu_m$. Observe that
we may suppose that $V_i$ do not intersect the divisor $D_m$. The
cyclic group $\mu_m$ acts on the cotangent sheaf of $D_{\zeta, \xi}$
thus on $\PP(\Omega^1_{D_{\zeta, \xi}})$. Denote by $\PP_{D_m}$ the
quotient $\PP(\Omega^1_{D_{\zeta, \xi}})/\mu_m$. There is a natural
projection $\PP_{D_m}\to V_i$ and the restriction of $\PP_{D_m}$ to
$V_i\setminus\{ R_i\}$ is isomorphic to the restriction of
$\Proj(\Omega_{X_m}^1(\log(D_m))$. Thus we can glue together the
restriction of $\Proj(\Omega_{X_m}^1(\log(D_m))$ to $X_m\setminus\{
R_i\}$ and  $\PP_{D_m}\to V_i$ to obtain a variety $\pi_m:\PP_m\to
X_m$ which is normal and $\Bbb Q$--factorial by construction
(locally it is quotient of a smooth variety by a finite group). One
easily verify that $\PP_m$ is projective and equipped with a $\Bbb
Q$--line bundle $\LL_m$ which has the searched properties. Observe
that $\PP_m$ has only isolated singular points.

The extension of the tautological inequality to singular variety is
straightforward: Let $X_m^{sm}$ be a desingularization of $X_m$,
$\PP_m^{sm}\to X_m^{sm}$ be the corresponding projective bundle of
the logarithmic differentials and $\LL_m^{sm}$ the tautological
bundle over it. Since $\PP_m$ and $\PP_m^{sm}$ are birational, there
exists a smooth variety $Z_m$, a commutative diagram

$$\matrix{ Z_m &\maprighto{a} &\PP_m^{sm}\cr
\mapdownr b & &\mapdownr{}\cr \PP_m &\maprighto{p_m}&B_m\cr}$$ where
the morphisms $a$ and $b$ are birational and a divisor $A$ on $Z_m$
such that $b^\*(\LL_m)=a^\*(\LL_M^{sm})+A$. Since the divisor $A$ is
vertical, (contained over a fibre of $p_m\circ b$) and the
tautological inequality holds on $X_m^{sm}$, the needed tautological
inequality holds on $X_m$.

The construction of the currents $T_m$ and $T'_m$ is similar to the
construction of the currents $T$ and $T'$: everything is defined to
let the construction work. Observe that the intersection of both
$T_m$ and $T'_m$ with a vertical divisor is zero.

To prove \ref{btildecurrent} we need the analogue of Stokes theorem
for currents:

\label stokesforcurrents. proposition\par\prop Let $T$ be a closed
positive current on a projective variety $X$. Let $f$ be a smooth
function on it. Then for almost all $\epsilon$ we there exists a
closed positive  current $T_\epsilon$ on $X_{\epsilon}:=\{ z\in
X\;\; / \;\; f(z)=\epsilon\}$ such that the following equality holds
for every smooth form $\omega$:
$$\int_{\{f\leq\epsilon\}}T\wedge
d(\omega)=\int_{X_\epsilon}T_\epsilon\wedge\omega.$$
\endprop

\noindent{\it Sketch of Proof}: We can find a sequence of smooth
closed currents $T_n$ such that $T_n\to T$ in the weak topology. By
Fubini theorem, we have that, for suitable $a$ and $b$ in $\Bbb R$
$$T(df\wedge d^c f)=\lim_{n\to\infty} T_n(df\wedge
d^cf)=\lim_{n\to\infty}\int_a^bdt\int_{X_t}T_n\wedge d^cf;$$ thus,
for almost all $\epsilon\in [a;b]$ the integrals
$\int_{X_t}T_n\wedge d^cf$ are uniformly bounded. Consequently, for
almost all $\epsilon\in[a;b]$ the measures $(T_n\wedge d^c
f)\vert_{X_\epsilon}$ on $X_\epsilon$ converge to a measure
$T_\epsilon\wedge d^cf$. If $\epsilon$ is outside the "bad set",the
classical Stokes theorem applied to the smooth closed currents gives
$$\int_{f\leq\epsilon}T\wedge d(\omega)=\lim_{n\to \infty}\int_{f\leq\epsilon}T_n\wedge d(\omega)
=\lim_{n\to \infty}\int_{X_\epsilon}T_n\wedge
\omega=\int_{X_\epsilon}T_\epsilon\wedge \omega.$$

\smallskip

Now we can give the

\noindent{\it Sketch of Proof of \ref{btildecurrent}}: Fix a smooth
metric on $\cO_X(B)$. Since $\II_BR=0$, by definition
$$R(B)=\lim_{\epsilon\to 0}\int_{\Vert
B\Vert\geq\epsilon}R\wedge c_1(\cO_X(D).$$ By Stokes theorem
\ref{stokesforcurrents}, for almost all $\epsilon$
$$\int_{\Vert
B\Vert\geq\epsilon}R\wedge c_1(\cO_X(D))=-\int_{\Vert
B\Vert\geq\epsilon}R\wedge dd^c\log\Vert B\Vert^2=\int_{\Vert
B\Vert=\epsilon}R_\epsilon\wedge{{d^c\Vert
B\Vert^2}\over{\epsilon^2}}\geq 0.$$ The conclusion follows.

\endssection

\endsection

\section Second approach to the theorems\par

\

In this section we will sketch the approach by Yamanoi to the main
theorems \ref{maitheoremalg} and \ref{maintheoreman}.

The Yamanoi approach is via the Ahlfors theory and works directly on
the moduli space of pointed stable curves of genus zero. The
complete proof requires a big machinery and is quite involved thus
we refer to the original paper [YA3] for the general statements. We
will give here a simplified proof, in the spirit of Yamanoi paper,
in the first non trivial case. The main ideas and difficulties
appear already here and we think that this case, and its proof, may
help to understand the general case.

The first step is the reduction to the case when $X$ is a blow up of
$\PP^1\times B$. This reduction goes back to Elkies [EL].

\label reduction. proposition\par\prop Suppose that
\ref{maitheoremalg} and \ref{maintheoreman} hold when $X$ is a blow
up of $\PP^1\times B$. Then \ref{maitheoremalg} and
\ref{maintheoreman} hold in general.
\endprop

\noindent{\it Sketch of Proof:} Let $(X;D)$ as in theorems
\ref{maitheoremalg} and \ref{maintheoreman}. changing $X$ by a
birational model of it, if necessary, we may suppose that there is a
generically finite morphism $g:X\to Z:=\PP^1\times B$ (commuting
with $p$) and a simple normal crossing divisor $H$ on $\PP^1\times
B$ such that $g^\*(K_Z(H))=K_X(D)+G$; where $G$ is a suitable
effective divisor on $X$ and (set theoretically) $g^{-1}(H)=D+G$.

Suppose that $f:Y\to X$ is a morphism from a curve, then since
\ref{maitheoremalg} or  \ref{maintheoreman} holds for $(Z,H)$, (we
omit $r$ in the analytic case) the inequality $(K_Z(H);Y)\leq
N^{(1)}_H(Y)+\chi(Y)+\epsilon (K_Z(H);Y)+\dots$ holds. Thus
$$(K_X(D);Y)+(G;Y)\leq
N^{(1)}_D(Y)+N^{(1)}_G(Y)+\chi(Y)+\epsilon(K_X(D+G);Y))+\dots.$$

We conclude because $N^{(1)}_G(Y)\leq (G;Y)$ and $K_X(D)$ is big.

\endsection

\section Ahlfors approach to SMT\par

\ssection Quick review of Ahlfors theory\par Suppose that $F$ and
$G$ are two bordered Riemann surfaces having finite Euler Poincar\'e
characteristic (eventually the boundary may be empty). Suppose that
$f:F\to G$ is an analytic finite morphism such that $f(\partial
F)\subset\partial G$ then the classical Hurwitz formula holds:
$$\chi(F)=\deg(f)\chi(G)+\sum_{P\in Ram(f)} (Ram_P(f)-1)$$ where
$Ram(f)$ is the set of ramification points of $f$ and $Ram_P(f)$ is
the ramification index of $f$ at $P$. Observe that we are using the
convention that $\chi(point)=-1$ or that $\chi(\PP^1)=-2$.

The first part of the Ahlfors theory is a generalization of this
formula when one removes the condition on the boundaries. Let $G^o$
be the interior of $G$; the set of points of the boundary of $F$
whose image is contained in $G^o$ is called the {\it relative
boundary of $f$}. Suppose that $H$ is a Riemann surface and $\eta$ a
pseudometric on it (i.e a smooth $(1,1)$ form which is positive
everywhere but a finite set of points where it vanishes); If $U$ is
a domain in $H$ we denote by $A(U,\eta)$ the area of $U$ with
respect to $\eta$; if $\beta$ is a Jordan curve on $H$ we denote by
$\ell(\beta,\eta)$ the length of $\beta$ with respect to the measure
defined by $\eta$; observe that they are both positive numbers.

We introduce a smooth positive metric $\omega$ on $G$ in such a way
that $A(G;\omega)<\infty$. The {\it mean sheet number of $f$} will
be the number $S_f:={{A(F;f^\*(\omega))}\over{A(G;\omega)}}$.
Observe that if $f$ is non ramified and unbordered, then $S_f$ is
the degree of $f$. If $U$ is a domain in $G$ then we define the {\it
sheet number of $U$ with respect to $f$}  by
$S_f(U):={{A(f^{-1}(U);f^\*(\omega))}\over{A(U;\omega)}}$.
Similarly, if $\beta$ is a Jordan curve on $G$, then we define the
{\it the sheet number of $\beta$} by
$L_f(\beta):={{\ell(f^{-1}(\beta);
f^{\ast}(\omega))}\over{\ell(\beta,\omega)}}$. We denote by $L_f$
the length of the relative boundary of $f$ with respect to
$f^\*(\omega)$. A morphism $f:F\to G$ will be said to be {\it
quasifinite} if it has finite fibres. The first main theorem of
Ahlfors theory is

\label ahlforscovering. theorem\par\thm Let $G$ be a bordered
Riemann surface, equipped with a positive metric $\omega$. Let $U$
be a domain and $\beta$ be a Jordan curve on $G$. Then there exist
positive constants $h$ and $k$ depending only on the metric and on
$U$ and $\beta$ respectively for which the following holds: For
every quasifinite morphism $f:F\to G$ from a bordered Riemann
surface to $G$ we have the following inequalities
$$\left\vert S_f-S_f(U)\right\vert \leq h L_f \;\; \; {\rm and
}\;\;\; \left\vert S_f-L_f(\beta)\right\vert \leq k L_f.$$
\endthm

For a proof we refer to [AH], to [HA] or to [NE].  What is very
important in the theorem above is that the constants $h$ and $k$
depend only on $U$ and $\beta$ (and on the metric $\omega)$ {\it but
not on } $F$ and $f$. The second main theorem of Ahlfors theory is
the following

\label ahlforsmaintheorem. theorem\par\thm Suppose that $G$ and $U$
is as in the previous theorem, then there is a constant $h>0$
depending only on $U$ (and the metric) such that, for every finite
covering $f:F\to G$ we have that
$$\max(\chi(f^{-1}(U));0)\geq \chi(U)S_f-hL_f.$$

\endthm

In the sequel we will denote by $a^+$ the number $\max(a,0)$.

\endssection

\smallskip

\ssection Ahlfors proof of SMT\par We will briefly show how to
deduce a form of the SMT from Ahlfors theorems. We will be a little
bit sketchy because these things are classical and well kown by
experts; we recall them here for reader's convenience and to point
out the analogies and the differences within the isotrivial and the
non isotrivial cases. Here and in the following we systematically
use the following:

\noindent -- We will always suppose that every (bordered) Riemann
surface we deal with will have finite Euler characteristic and it is
either compact or  it is relatively compact in a bigger Riemann
surface.

\noindent -- {\it Mayer--Vietoris formula}: If F is a Riemann
surface and $U$ and $V$ are two open sets of $F$ then
$\chi(F)=\chi(U)+\chi(V)-\chi(U\cap V)$.

\noindent -- If $\beta$ is a non compact Jordan curve which divides
$F$ in two connected components $U$ and $V$ then
$\chi(F)=\chi(U)+\chi(V)+1$. We will call $\beta$ a {\it cross cut}.

\noindent -- The Euler--Poincar\'e characteristic of a connected
Riemann surface is at least $-2$ and it is $-2$ if and only if it is
isomorphic to $\PP^1$.

\noindent -- Let $f:F\to G$ be a finite covering, Let $U$ be a
domain in $G$. A connected component $V$ of $f^{-1}(U)$ is called a
{\it island} if it is relatively compact in $F$ and a {\it
peninsula} otherwise.

\smallskip

Suppose that $P_1,\dots, P_q$ are $q$ points on $\PP^1$ and $U_1,
\dots, U_q$ are small disks around the $P_i$'s whose the closure are
mutually disjoint. Denote by $G^0$ the Riemann surface
$\PP^1\setminus\bigcup_iU_i$. We fix on $\PP^1$ the Fubini--Study
metric $\omega_{FS}$: $A(\PP^1;\omega_{FS})=1$.

If $f:F\to \PP^1$ is a  quasifinite morphism, then we denote by
$N_i(f)$ the number of islands on $F$ above $U_i$. The theorem which
generalize the SMT is the following, it can be seen as a strong, non
integrated form of it.

\label ahlforssmt. theorem\par\thm Suppose that we fixed $U_i$  as
above, then there is a positive constant $h$ depending only on the
$U_i$'s such that the following holds: for every Riemann surface $F$
and quasifinite morphism $f:F\to \PP^1$ we have that
$$\chi^+(F)+\sum_iN_i(f)\geq (q-2)A(F;f^\*(\omega_{FS}))-hL_f.$$
\endthm

Theorem \ref{ahlforssmt} is a consequence of \ref{ahlforscovering}
and \ref{ahlforsmaintheorem}. We give here a Sketch of the proof;

\noindent{\it Sketch of Proof:} Denote by $G_0$ the open set
$\PP^1\setminus\bigcup_{i=1}^q\overline{U}_i$ and $\beta$ the
boundary of $G_0$. The Euler characteristic of $G_0$ is $q-2$.
Denote by ${\cal I}$ (resp. ${\cal P}$) the set of islands (resp.
peninsulas) of $F$ over the $U_i$. Let $F_0$ be $f^{-1}(G_0)$ and
$\gamma=f^{-1}(\beta)$. By Mayer Vietoris Formula, we have
$$\chi(F)=\chi(F_0)+\sum_{I\in{\cal I}}\chi(I)+\sum_{P\in{\cal
P}}\chi(P)+n;$$ where $n$ is the number of cross cuts of $\gamma$
(components which are not compact). Since for every connected
component $A$ in the sum, $\chi(A)\geq-1$, each peninsula touch at
least a cross cut and each cross cut touch at most one peninsula,
$$\chi^+(F)+\sum_iN_i(f)\geq\chi^+(G_0).$$
We conclude applying \ref{ahlforscovering} and
\ref{ahlforsmaintheorem}.

Denote by $n(f,P_i)$ the cardinality of the $z\in F$ such that
$f(z)=P_i$ then one easily sees that
$\sum_in(f,P_i)\geq\sum_iN_i(f)$.

Let $(Y,g)$ be a parabolic Riemann surface and $f:Y\to \PP^1$ an
analytic map. Apply the theorem to $F_t:=\{ z\in Y \; {\rm s.t.}\;
g(z)\leq t\}$. It is well known that
$$\lim_{t\to\infty}{\int_1^r{{L_{f_t}dt}\over{t}}\over{\int_1^r{{A(F_t,f^\ast(\omega_{FS}))dt}\over{t}}}}=0;$$
where $L_{f_t}$ is the length of the relative boundary of $F_t$ (cf.
for instance [BR]). Thus if one integrate the inequality of the
theorem with respect to $\int_1^r{{dt}\over{t}}$, one finds a proof
of the SMT.

Remark that in the proof we are allowed to move a little bit the
points $P_i$'s and the results remains unchanged! This means that
the SMT remains true if we perturb a little bit the divisor
$D:=\sum_iP_i$. This is the key point of the Yamanoi approach: In
the theorem we can move a little bit the divisor and everything
remains true, thus we can give a general proof working on the moduli
space of stable pointed curves of genus zero, which is compact! The
only problem is that sometimes the points $P_i$'s may coincide.

One works directly on the moduli space of stable curves of genus
zero with $n$ marked points $\M_{0,n}$ and on its universal family
$p:\U_{0,n}\to \M_{0,n}$. It is well known that there are $n$
sections $\xi_i:\M_{0,n}\to \U_{0,n}$ of $p$ and that ${\cal
D}_n:=\sum_i\xi_i(\M_{0,n})$ is the universal divisor: the
restriction of ${\cal D}_n$ to the generic fibre of $p$ is the
divisor given by the marked points. Let $K_{\U/\M}$ be the relative
dualizing sheaf of $p$. In the sequel we will denote by $K_n$ the
line bundle $K_{\U/\M}({\cal D}_n)$ on $\U_{0,n}$; we will suppose
that it is equipped with a smooth hermitian metric and we will
denote by $\omega$ its first Chern form.

Since a rigorous proof is quite involved and requires a careful
attention to details, we will explain the main steps of the proof in
the case when $n=4$ (stable curves of genus zero with $4$ marked
points). We refer to the original paper by Yamanoi for the general
case. This is the first non trivial case which cannot be deduced
directly from the classical SMT. Even in this case a detailed proof
requires a skillful work (we think that filling the gaps is a good
exercise). Nevertheless we think that all the main steps and ideas
of the proof are already present in this case.

\smallskip

\ssection Explicit description of $\M_{0,4}$ and $\U_{0,4}$\par The
moduli space $\M_{0,4}$ is isomorphic to the projective line
$\PP^1$.

Let $X:=\PP^1\times\PP^1$; we denote by $p:X\to\PP^1$ the first
projection. The map $p$ is equipped with $4$ sections: we write them
in affine coordinate: $\xi_0(z):=(z,0)$, $\xi_1(z)=(z,1)$,
$\xi_{\infty}(z)=(z,\infty)$ and $\xi_\Delta(z)=(z,z)$; we will
denote by $\xi_i$ and $\xi_\Delta$ the image of the $\xi_i$ and of
$\xi_\Delta$ respectively. The $\xi_i$, for $i=0, 1,\infty$, do not
intersect and  $\xi_\Delta$ intersect the  $\xi_i$ properly over
$i$.

Let $\pi:\tilde X\to X$ be the blow up of $X$ over the three points
$\xi_\Delta\cap\xi_i$. Then $\tilde X$ is the universal family
$\U_{0,4}$ and the strict transforms $\hat\xi_j$'s of the $\xi_j$'s,
for $j=0,1,\infty$ and $\Delta$ are the universal sections. The map
$p\circ\pi: \tilde X\to\PP^1$ is the universal map.

Let $U_g\subset\PP^1$ be the open set $\PP^1\setminus\{
0,1,\infty\}$. Then $\tilde X\vert_{U_g}:=p^{-1}(U_g)$ is isomorphic
to $U_g\times\PP^1$; let $h:\tilde X\vert_{U_g}\to\PP^1$ be the
second projection. Suppose that $z=0,1$ or $z=\infty$, then we can
find a neighborhood $U_x\subset\PP^1$ of $z$ for which $\tilde
X\vert_{U_z}:=p^{-1}(U_z)$ is the blow up of $U_z\times\PP^1$ over
the point $(z,z)$; by construction there is a projection
$h:(g_1;g_2):\tilde X\vert_{U_x}\to\PP^1\times\PP^1$. It is easy to
see (by restriction to the fibres of $p$) that the restriction of
$K_4$ to $\tilde X\vert_{U_g}$ is $h^\*(\cO(2))$. Moreover for
$x=0,1,\infty$,  (we may suppose that) the restriction of $K_4$ to
$\tilde X\vert_{U_z}$ is $h^\*(\cO(1,1))$. In the sequel we will
suppose that the metric on the restriction of $K_4$ to these open
sets is the pull back of the Fubini--Study metrics; this is not
exactly the case but since $\tilde X$ is compact, the error we make
is bounded and can be controlled.

\endssection

In the sequel we will suppose that we are in the following
situation: $R$ will be a open set in $\PP^1$ (for the analytic
topology). $g:F\to R$ is a proper maps between Riemann surfaces and
$f:F\to \tilde X$ is an analytic map such that the following diagram
is commutative:
$$\matrix{F & \llongmaprighto{f} &\tilde X\cr \mapdownl{g}& &
\mapdownr{\pi}\cr R &\llongmaprighto{\iota}&\PP^1;\cr}$$ where
$\iota:R\to\PP^1$ is the inclusion. We will call this {\it a
situation}.

Suppose we are in a situation as above, and $W\subseteq F$ is a open
set, we will denote by $A(W,\omega)$ the area of $W$ with respect to
the volume form $f^{\ast}(\omega)$ on $F$. If $\gamma$ is a Jordan
curve on $F$, we denote by $L(\gamma, \omega)$ the length of
$\gamma$ with respect to the measure defined by $f^{\ast}(\omega)$.
We will denote by $L_f$ the length of the relative boundary of $f$.

Let $D={\cal D}_4\hookrightarrow \tilde X$ be the universal divisor:
we will denote by $n(D,f)$ the cardinality of the set $\{ z\in F\;
/\; f(z)\in D\}$.

\endsection

\

\section The local version of the theorem\par

\

The key step of Yamanoi proof is a local version of the theorem.
This local version plays the role of theorem \ref{ahlforssmt} in the
Ahlfors proof of SMT. Given $f:F\to \tilde X$,  and a open set $U$
of $\tilde X$, we will generalize in the obvious way the concept of
island and peninsula of $F$ over $U$: an island will be a connected
component of $f^{-1}(U)$ which is relatively compact, etc. We will
denote by $N(f,U)$ the number of islands of $F$ over $U$.

Before we state and prove the theorem, we need to state a
generalization of a classical theorem by Rouch\'e:

\label rouche. proposition\par\prop Let $E$ be a Jordan domain of
$\PP^1$ and $b\in E$; then there exists a positive constant
$C:=C(E,b)$ with the following property: Let $F$ be a bordered
Riemann surface and $\zeta:F\to E$ an analytic function such that
$\zeta(F)=E$ and $\zeta(\partial(F))=\partial E$; then for every
meromorphic function $\alpha:F\to\PP^1$ such that $\vert
\alpha(z)-b\vert<C$ for every $z\in F$, there exists $z\in F$ with
$\alpha(z)=\zeta(z)$.
\endprop

The proof of this proposition is a variation of the classical
Rouch\'e theorem and can be found on [YA2].

The local version of Yamanoi theorem is

\label localyamanoi. theorem\par\thm Let $x\in\PP^1$ then we can
find a open neighborhood $x\in U_x\subseteq\PP^1$,  open
neighborhoods $W_i\subseteq \tilde X\vert_{U_x}$ of $\xi_i\cap\tilde
X\vert_{U_x}$, for $i=0,1,\infty$ and $\Delta$, with disjoint
closures, and a positive constant $h_x$ for which the following
holds:

For every situation
$$\matrix{F & \llongmaprighto{f} &\tilde X\cr \mapdownl{g}& &
\mapdownr{\pi}\cr R &\llongmaprighto{\iota}&\PP^1;\cr}$$  for which
$x\in R\subseteq U_x$,  we have that
$$h_xL_f+\deg(g)+\chi^+(F)+\sum_iN(f,W_i)\geq A(F,\omega).$$
\endthm

We recall that $L_f$ is the length of the relative boundary. One
sees the similarity of the theorem above with the classical theorem
by Ahlfors \ref{ahlforssmt}. One should notice that \ref{ahlforssmt}
is one of the main tools of the proof of the theorem above.

\Proof First case: we suppose that $x\neq 0,1$ or $\infty$. In this
case the theorem is essentially \ref{ahlforssmt}; we give some
details: Take a small disk $U_x$ around $x$; then $\tilde
X\vert_{U_x}$ is isomorphic to $U_x\times\PP^1$; let $h:\tilde
X\vert_{U_x}\to\PP^1$ be the second projection. We may suppose that,
for $i=0,1$ and $\infty$, we have $h\circ\xi_i(x)=i$ and
$\xi_{\Delta}(x)=x$. Take small neighborhoods $U_i$ of $i$ in
$\PP^1$ and a small neighborhood $U_\Delta$ of $x$ with non
intersecting closures. We obtain the theorem in this case by
applying Ahlfors theorem \ref{ahlforssmt} to the morphism $h\circ
f:F\to \PP^1$. Notice that in this case the term $\deg(g)$ is not
there.

The new case is when $x=0,1$ or $\infty$.

Suppose that $x=0,1$ or $\infty$: we may suppose that $x=0$ the two
other cases are similar.

In this case the fibre of $\pi:\tilde X\to \PP^1$ over $x$ is the
union of two components $X_1$ and $X_2$ both isomorphic to $\PP^1$.
The universal sections $\xi_0$ and $\xi_\Delta$ intersect $X_1$ and
not $X_2$ while $\xi_1$ and $\xi_\infty$ intersect $X_2$ and not
$X_1$. Take a neighborhood $U_x$ of $x$ and two maps $h_i:\tilde
X\vert_{U_x}\to \PP^1$. We may suppose that: $h_1(\xi_\infty(x))=0$,
$h_1(\xi_1(x))=1$, $h_1(X_2)=\infty$ and that $h_2(\xi_0(x))=0$,
$h_2(\xi_\Delta(x))=1$, $h_2(X_1)=\infty$.

For $i=0,1$ and $\infty$, choose small neighborhoods $U_i$ of $i$
whose closure do not intersect. Call $A:=h_1^{-1}(U_\infty)$ and
$B:=h_2^{-1}(U_\infty)$. We may take the neighborhood $U_x$ so small
that $h_1(A\cap B)=h_2(A\cap B)=U_{\infty}$.  We may also suppose
that $U_x$ is so small that $\vert
h_j(\xi_i(z))-\ell\vert<C(i,U_\ell)$  for $z\in U_x$, $j=1$, $i=0,
1, \Delta$ $\ell= 0,\infty, 1$ respectively, or $j=2$ and $i=\ell=
0,1,\infty$ respectively  (cf. prop. \ref{rouche}).

For $i=0,1$ define the following open sets $W_i:=h_1^{-1}(U_i)$ and
$V_i:=h_2^{-1}(U_i)$. Notice that the $W_i$ and $V_j$ are mutually
disjoint and we may suppose $U_x$ so small that
$\xi_\infty(U_x))\subset W_0$, $\xi_1(U_x))\subset W_1$,
$\xi_0(U_x)\subset V_0$ and $\xi_\Delta(U_x)\subset V_1$.

Let $A_1,\dots, A_r$ be the island of $F$ over $A\cap B$,
$F_1:=F\setminus\{A_1,\dots, A_r\}$  and ${\cal I}$ and ${\cal P} $
be the set of islands and peninsulas of $F$ over $B$ respectively.
Remark that ${\cal I}$ is also the set of islands of $F_1$ over
$U_{\infty}$ via $h_1\circ f$.

Let $S$ be one of the $A_i$'s or an element of one of ${\cal I}$.
Prop. \ref{rouche} applied to $S$, $h_1\circ f\vert_{S}$ and
$h_1\circ\xi_1\circ g$ implies that there is $z\in S$ such that
$h_1\circ f(z)=h_1\circ\xi_1\circ g(z)$. Since $h_1$ restricted to a
fibre of $\pi$ different from the fibre over $x$ is an isomorphism,
one finds that $\pi(f(z))=x$. Consequently every such island
intersects the fibre over $x$ and the properness of $g$ implies that
$r+Card({\cal I})\leq\deg(g)$. In particular \labelf
chiofF1\par$$\chi(F_1)+Card({\cal I})
\leq\chi(F)+\deg(g).\eqno{{(\numfo)}}$$\advance\ssnu by1

Let $N_A$ (resp. $N_B$) the number of islands of $F_1$ (or $F$ which
is the same) over $V_0$ and $V_1$(resp. over $W_1$ and $W_2$).  A
direct application of Ahlfors theory and Mayer--Vietoris formula to
$F_1$ and $h_2\circ f$ gives a universal constant $h$ (independent
on $F$) such that \labelf
yamalhfors1\par$$\chi^+(F_1)+N_A-\sum_{P\in{\cal
P}}\chi^+(P)-\sum_{I\in{\cal I}}\chi(I)\geq A(F;(h_2\circ
f)^{\ast}(\omega_{FS}))-hL_f.\eqno{{(\numfo)}}$$ Here and in the
sequel, we systematically use theorem \ref{ahlforscovering}.

We apply again Ahlfors theory to each island and peninsula of $F_1$
over $B$. Observe that, for each island $I$, $\chi^+(I)<\chi(I)+1$.
Thus we obtain \advance\ssnu by1\labelf
yamahlfors2\par$$\sum_{P\in{\cal P}}\chi^+(P)+\sum_{I\in{\cal
I}}\chi(I)+Card({\cal I})+ N_B\geq A(F;(h_1\circ
f)^{\ast}(\omega_{FS}))-hL_f.\eqno{{(\numfo)}}$$

The conclusion follows from \ref{chiofF1}, \ref{yamalhfors1} and
\ref{yamahlfors2}.

Remark that, since the base $\PP^1$ is compact, the error we make
using the pull back  of the Fubini Study metric via $h_i$ instead of
the $(1,1)$ form $\omega$ of $K_4$ over $\tilde X$ is controlled by
changing the constant $h_x$.

If we have a {\it situation} as above, we will denote by $n(D,f)$
the number of points $z\in F$ such that $f(z)\in D$. Let $R_f$ be
the number $\sum_{z\in F}(Ram(g)-1)$

As a consequence, we find

\label localyamanoi2. theorem\par\thm Suppose that the we are in the
hypotheses of theorem \ref{localyamanoi}. Then
$$A(F,\omega)\leq n(D,f)+ R_f+ \deg(g)\chi(R)+\deg(g)+h_xL_F.$$
\endthm

\Proof Since $g$ is proper, $\chi^+(F)\leq \chi(F)+\deg(g)$ and by
Hurwitz formula, $\chi(F)=\deg(g)\chi(R)+R_f$. Thus it suffices to
apply \ref{localyamanoi} and \ref{rouche} to $f$ and $\xi_i\circ g$
over each island.

\smallskip

Observe that the theorem above is a local version of the theorems;
It seems better then the theorem because one has the impression that
one can put $\epsilon=0$; nevertheless there is the term coming from
the relative boundary $L_f$. We will see in the sequel we will need
to put $\epsilon>0$ in order to control this term. Even if this is
not the only reason, it is the most important.

\endssection

\endsection

\

\section The non integrated version of the theorem\par

\

After the local version of the theorem we will prove a global non
integrated version of the theorem. Here too we will put some
restrictive hypotheses on the situations: nevertheless we would like
to remark that these hypotheses are suffice to prove theorems
\ref{maitheoremalg} and \ref{maintheoreman}.

Let $K$ be a compact set of $\PP^1$ (which may be empty) We will say
that a sequence of open sets $R_n\subset \PP^1$ is {\it relatively
exhausting with respect to $K$} if, for every compact disk
$\Delta\subset\PP^1\setminus K$, there exists $n_0$ such that, for
every $n\geq n_0$ we have that $\Delta\subset R_n$.

The non integrated version of the theorems we propose is the
following

\label nonintyama. theorem\par\thm Let $K$ be a compact set of
$\PP^1$ and $\epsilon>0$. Suppose that
$$\matrix{F_n & \llongmaprighto{f_n} &\tilde X\cr \mapdownl{g_n}& &
\mapdownr{\pi}\cr R_n &\llongmaprighto{\iota_n}&\PP^1;\cr}$$ is a
sequence of situations with the sequence $\{ R_n\}$ relatively
exhausting with respect to $K$. Then, after subsequencing, we can
find  constants $h$ and $C$ such that, for every term of the
subsequence
$$A(F_n;\omega)\leq n(D, f_n)+R_{f_n}+\epsilon A(F_n,\omega)+h\ell(\partial
F_n,\omega)+\deg(g_n)(\chi(R_n)+C).$$ The constant $h$ is
independent on the sequence and $C$ depends only on the sequence
(and not on the terms of the sequence).
\endthm

\Proof We can find a open set $W$ containing $K$ having the
following property: the open set $W\setminus K$ is a finite union of
open sets of the form $U_x$ of theorem \ref{localyamanoi} such that
$U_x\cap K\neq\emptyset$ and $x\not\in K$.  Choose an integer
$J>{{2}\over{\epsilon}}$. Let $\gamma_1,\dots, \gamma_J$ be Jordan
curves of $\PP^1$  and $\delta_i$ small open neighborhoods of
$\gamma_i$ for which the following properties hold:

\noindent -- For each $i$, every connected component of
$\PP^1\setminus\gamma_i$ is simply connected and contained in one of
the open sets $U_x$ of theorem \ref{localyamanoi}.

\noindent -- For each $i$, every connected component of
$\PP^1\setminus\delta_i$ is again simply connected (and contained in
one of the $U_x$).

\noindent -- For every triple of distinct indices $(i ,j ,k)$ we
have $\delta_i\cap\delta_j\cap\delta_k=\emptyset$.

\noindent -- If a connected component of $\PP^1\setminus\delta_i$
intersects $W$ then it is contained in it.

Because of the third condition, for every $n$, we have
$$\sum_jA(g_n^{-1}(\delta_j);\omega)\leq 2A(F_n,\omega).$$
Thus we can find a $j_0$ and a subsequence $n_k$ for which \labelf
areaofdelta\par$$A(g_{n_k}^{-1}(\delta_{j_0}),\omega)\leq
{{2}\over{J}}A(F_{n_k};\omega)\leq\epsilon
A(F,\omega).\eqno{(\numfo)}$$ Fix such a $j_0$ and call $\delta$ the
open set $\delta_{j_0}$ etc. From now on, we will omit to change
notation when we pass to a subsequence.

Let $U$ be a connected component of $\PP^1\setminus\delta$ and
consider the set of Riemann surfaces $F_{n,U}:=g_n^{-1}(R_n\cap U)$.

Either $\limsup_n{{A(F_{n,U},\omega)}\over{\deg(g_n)}}<\infty$ or
$\limsup_n{{A(F_{n,U},\omega)}\over{\deg(g_n)}}=\infty$. We suppose
that we are in the second situation, thus, passing to a subsequence,
we may suppose that the $\limsup$ is indeed a limit.

If $U$ is not contained in $W$ then we may also suppose that $R_n$
contains $U$ for every $n$. Suppose that we are in this case.

Denote by $\Delta_r$ the disk of radius $r$. We may suppose that $U$
is biholomorphic to the disk $\Delta_{r_0}$ for some $r_0<1$. We may
also suppose that $U\simeq \Delta_{r_0}\subset\Delta_1\subseteq U_x$
for some $x\in U$. Let $F_{n,\Delta}:=g_n^{-1}(\Delta_1)$ and for
every $r\in (0,1)$ let $F_{n,r}:=g_n^{-1}(\Delta_r)\subset
F_{n,\Delta}$ . We can find a non negative function $G$ which is
$C^{\infty}$ outside the ramification points of $g_n$ and integrable
on $F_{n,\Delta}$, such that
$f_n^{\ast}(\omega)\vert_{F_{n,\Delta}}=\sqrt{-1}G^2dg_n\wedge
d\overline{g}_n$. Let $$S_n(r):=\int_{0}^rdt\int_{\partial
F_{n,t}}Gtd\arg(g_n),$$ then ${{dS_n}\over{dr}}=\ell(\partial
F_{n,r},\omega)$. By Cauchy Schwartz inequality we have
$$\eqalign{S_n(r)&\leq \left(\int_0^rtd\wedge
d\arg{g_n}\right)^{1/2}\cdot\left(\int_0^rG^2tdt\wedge
d\arg(g_n)\right)^{1/2}\cr &=C_U\cdot\left(
\deg(g_n)\right)^{1/2}\cdot\left(
A(F_{n,r};\omega)\right)^{1/2}.\cr}$$ Where $C_U>0$ is a constant
depending only on $U$.

\label measurelemma. lemma\par\lemma Let $\delta>0$ and $S_n(r)$ be
a sequence of differentiable functions on $[0,1)$ with $S_n(r_0)\geq
n^{2/\delta}$ then the set
$$I_S:=\left\{ 1>r\geq r_0 \; /\; S'_n(r)\geq
{{S_n(r)}\over{(1-r^2)}}\; {\rm for\; some\; } n\right\}$$ is such
that $\int_{I_S}{{dr}\over{1-r^2}}<\infty$.

\endlemma

The proof of the lemma is standard and can be found on [MQ2].

As a consequence of the lemma above, for every $\epsilon'>0$ we can
find a subsequence of the $F_n$ and a $R>r_0$ for which
$\ell(\partial F_{n,R},\omega)<\epsilon' A(F_{n,R}, \omega)$. We
call again $U$ (by abuse of notations) the enlarged open set for
which this last inequality holds.

A similar argument holds when $U$ is contained in $W$. In this case,
even taking a subsequence, we cannot suppose that $U$ is contained
in $R_n$:  enlarging a bit $U$, as before, we may suppose
$$\ell(\partial F_{n,U}, \omega)\leq \epsilon
A(F_{n,U},\omega)+\ell(\partial F_n\cap F_{n,U}; \omega).$$

Since $F_n=\bigcup_Ug_n^{-1}(U)\cap g_n^{-1}(\delta)$, and
\ref{areaofdelta} holds, we apply theorem \ref{localyamanoi2} and
obtain
$$\eqalign{&A(F_n;\omega)\leq\sum_{U}A(F_{n,U},\omega)+A(g_n^{-1}(\delta))\cr
& \leq \sum_{U}\left(n(D,f_n\vert_{F_{n,U}})+
R_{f_n\vert_{F_{n,U}}}+\deg(g_n)\chi(R_n\cap U)+\deg(g_n)+\epsilon
A(F_{n,U},\omega)\right)\cr &+\epsilon A(F_n,\omega)+ h\ell(\partial
F_n,\omega)+ C,\cr}$$ Where the constant $C$  take care of the open
sets $U$ for which
$\limsup_n{{A(F_{n,U},\omega)}\over{\deg(g_n)}}<\infty$ thus depends
only on the sequence.

Since $\chi(R_n)\geq \sum_U\chi(R_n\cap U)$ we conclude that
$$A(F_n;\omega)\leq n(D, f_n)+R_{f_n}+\epsilon A(F_n,\omega)+ h\ell(\partial
F_n,\omega)+\deg(g_n)(\chi(R_n)+C).$$

From this we deduce

\label nonintyama2. theorem\par\thm Let $\epsilon>0$ then there
exists constants $C$ and $h$ such that, {\rm for every} situation as
above,
$$A(F;\omega)\leq n(D, f)+R_{f}+\epsilon A(F,\omega)+ h\ell(\partial
F,\omega)+\deg(g)(\chi(R)+C).$$
\endthm

\Proof If not, we can find a sequence of situations for which
$$\lim_{n\to\infty}{{1}\over{\deg(g_n)}}\cdot\left((1-\epsilon) A(F_n;\omega)- \left(n(D, f_n)+R_{f_n}+ h+ h\ell(\partial
F_n,\omega)\right)\right)+\chi(R_n)=+\infty.$$ And this contradicts
theorem \ref{nonintyama}.

This easily imply, together with \ref{reduction} the algebraic
version of $abc$. The analytic version of $abc$ requires again a
control of the length of the boundary; this is again standard: We
give a sketch of the proof in  a special case. We suppose, to
simplify that $g:Y\to\CC$ is a proper map and the following diagram
is commutative

$$\matrix{Y & \llongmaprighto{f} &\tilde X\cr \mapdownl{g}& &
\mapdownr{\pi}\cr \CC&\llongmaprighto{\iota}&\PP^1.\cr}$$

Applying theorem \ref{nonintyama2} when  $R=R_t$, the disk of radius
$t$ in $\CC$,  and integrating  with respect to
$\int_0^r{{dt}\over{t}}$, we obtain
$$(K_4, Y)(r)\leq N_{D_4}^{(1)}(Y)(r) +R_f(r)+\epsilon(K_4,Y)(r) +
h\int_{0}^r{{\ell(\partial g^{-1}(R_t),\omega)dt}\over{t}}+ C\log
r.$$

We can write $f^{\ast}(\omega)=\sqrt{-1}G^2dg\wedge d\overline{g}$
with $G$ a non negative function which is integrable and
$C^{\infty}$ outside the ramification points of $g$.

Introduce the function
$S(r):=\int_0^{\log(r)}dt\int_{g=t}tGd\arg(g)$. We have that
$S'(r)=\int_{0}^r{{\ell(\partial g^{-1}(R_t),\omega)dt}\over{t}}$
and Cauchy--Schwartz inequality gives $S(r)\leq
C(\log(r))^{1/2}\cdot {{d(K_4;Y)(r)}\over{dr}}$. A double
application of  lemma \ref{lang2} allows to conclude the proof.
Remark that this argument is similar to the argument used to derive
the SMT from Ahlfors theory.

\smallskip

\ssection The general statement proved by Yamanoi\par As a
conclusion, we state without proof the main theorem proved in [YA3];
we refer to the original paper for the proof. Let $n>3$ be an
integer. A {\it situation} will be a commutative diagram
$$\matrix{F & \llongmaprighto{f} &{\cal U}_{0,n}\cr \mapdownl{g}& &
\mapdownr{\pi}\cr R &\llongmaprighto{\iota}&{\cal M}_{0,n};\cr}$$
with $F$ and $R$ bordered Riemann surfaces, $g$ a proper analytic
map and $f$ and $\iota$ analytic. We fix a metric on $K_n$  and a
positive $(1,1)$ form $\eta$ on ${\cal M}_{0,n}$. We define $n({\cal
D}_n,f)$ and $R_g$ as before. Observe that ${\cal M}_{0,n}={\cal
U}_{0,n-1}$ thus we may define $n({\cal D}_{n-1},\iota)$.

\label generallocalyama. theorem\par\thm Let $\epsilon>0$ then there
is a constant $C$ depending only on $\epsilon$ and the metrics
chosen on $K_n$ and ${\cal M}_{0,n}$ with the following property:
For every situation as above, we have
$$\eqalign{A(F, f^{\*}&(c_1(K_n)))\leq n({\cal D}_n,f)+ R_g+\epsilon
A(F,f^{\*}(c_1(K_n)))\cr &+ C\deg(g)\left( A(R,g^{\*}(\eta))+n({\cal
D}_{n-1},\iota)+\chi^+(R)+\ell(\partial F,
f^{\*}(c_1(K_n)))\right).\cr}$$

\endthm
Cf. [YA3] Theorem 4. An argument similar to the one sketched above
allows to deduce the $abc$ conjecture from the theorem above.

\endssection

\endsection

\

\section Conclusions and final observations\par

\

A posteriori one would like to compare the two proofs. The proof by
McQuillan has a global nature while the Yamanoi 's is more local. Of
course one is tempted to apply the techniques to other situations;
for instance to families of surfaces of general type. The first part
of the proof by McQuillan passes through without pain (essentially
everything until prop. \ref{pminusdelta}). Then one have to deal
with a more subtle situation: here we strongly used the fact that
the singularities of families of semistable curves are well
understood and quite easy. In general the situation is more
complicated.

The Yamanoi approach is essentially local. Suppose that we have a
family of varieties over a curve and we want an inequality similar
to the $abc$ in this situation. Split the base in finitely many
small open sets $U_i$. Take a sequence of curves with maps in our
family. Look to the sequence of the areas of the preimages of each
$U_i$. If this is bounded, there is nothing to prove. If it is
unbounded, then one look for a local inequality which will involve
the length of the boundary as in theorem \ref{localyamanoi2}. Then
one can conclude adapting the arguments of theorem \ref{nonintyama}.
Of course this will need a generalization of Ahlfors theory (even in
the smooth case) and to our knowledge this is still unknown.

In conclusion, the first part of the proof by McQuillan and the
second part of the proof by Yamanoi can be generalized.  Thus it is
probable that the best way to proceed will be by applying a mix of
both proofs!
\endsection

\

\

\centerline {\sectionfont References.}

\

\spacing \item{[AH]} Ahlfors, Lars, Zur Theorie der
\"Uberlagerungsfl\"achen. Acta Math. 65 (1935), no. 1, 157--194.

\spacing \item{[AS]} Ahlfors, Lars V.; Sario, Leo, Riemann surfaces.
Princeton Mathematical Series, No. 26 Princeton University Press,
Princeton, N.J. 1960 xi+382 pp.

\spacing \item{[BG]} Bombieri, Enrico; Gubler, Walter, Heights in
Diophantine geometry. New Mathematical Monographs, 4. Cambridge
University Press, Cambridge, 2006. xvi+652 pp.

\spacing \item{[BR]} Brunella, Marco, Courbes entières et
feuilletages holomorphes. Enseign. Math. (2) 45 (1999), no. 1-2,
195--216.

\spacing\item{[CH]} Chen, Xi, On Vojta $1+\epsilon$ conjecture.
prepint avaible at arXiv:0705.1727.

\spacing \item{[EL]} Elkies, Noam D. $ABC$ implies Mordell.
Internat. Math. Res. Notices 1991, no. 7, 99--109.

\spacing\item{[FA]} Faltings, G. Endlichkeitss\"atze f\"ur abelsche
Variet\"aten ber Zahlk\"orpern. Inv\-ent. Ma\-th. 73 (1983), no. 3,
349-366.

\spacing \item{[GK]} Griffiths, Phillip; King, James Nevanlinna,
theory and holomorphic mappings between algebraic varieties. Acta
Math. 130 (1973), 145--220.

\spacing \item{[HA]} Hayman, W. K. Meromorphic functions. Oxford
Mathematical Monographs Clarendon Press, Oxford 1964 xiv+191.

\spacing\item{[HS]} Hindry, M.; Silverman, J. H. The canonical
height and integral points on elliptic curves. Invent. Math. 93
(1988), no. 2, 419--450.

\spacing \item{[KI]} Kim, Minhyong, Geometric height inequalities
and the Kodaira-Spencer map. Compositio Math. 105 (1997), no. 1,
43--54.

\spacing \item{[MQ1]} McQuillan, Michael, Diophantine approximations
and foliations. Inst. Hautes Études Sci. Publ. Math. No. 87 (1998),
121--174.

\spacing \item{[MQ2]} McQuillan, Michael, Non commutative Mori
theory, Preprint IHES

\spacing \item{[MQ3]} McQuillan, Michael, Old and new techniques in
function fields arithmetics, preprint

\spacing \item{[MQ4]} McQuillan, Michael, Rational criteria for
hyperbolicity. book preprint.

\spacing\item{[MY]}  Miyaoka, Yoichi, The orbibundle
Miyaoka-Yau-Sakai Inequality and an effective
Bo\-go\-mo\-lov--McQuillan Theorem, preprint.

\spacing \item{[NE]}Nevanlinna, Rolf, Analytic functions. Die
Grundlehren der mathematischen Wissenschaften, Band 162
Springer-Verlag, New York-Berlin 1970 viii+373 pp.

\spacing \item{[NI]} Nitaj, Abderrahmane, the $abc$ conjecture web
page:

http://www.math.unicaen.fr/$\;\tilde{ }$ nitaj/abc.html

\spacing \item{[NWY]} Noguchi, Junjiro; Winkelmann, J\"org; Yamanoi,
Katsutoshi, The second main theorem for holomorphic curves into
semi-abelian varieties. Acta Math. 188 (2002), no. 1, 129--161.

\spacing\item{[OE]} Oesterl\'e, Joseph, Nouvelles approches du
"th\'eorème" de Fermat.  S\'eminaire Bourbaki, Vol. 1987/88.
Ast\'erisque No. 161-162 (1988), Exp. No. 694, 4, 165--186 (1989).

\spacing \item{[OS]} Osgood, Charles F. Sometimes effective
Thue-Siegel-Roth-Schmidt-Nevanlinna bou\-nds, or better. J. Number
Theory 21 (1985), no. 3, 347--389

\spacing \item{[SA]} Sauer, Andreas, Deficient rational functions
and Ahlfors's theory of covering surfaces. Ark. Mat. 39 (2001), no.
1, 151--155.

\item{[SE]} Serre, J.P. Lectures on Mordell Weil Theorem.
Friedr. Vieweg \& Sohn, (1989).

\spacing\item{[SZ]}  Szpiro, L. Discriminant et conducteur des
courbes elliptiques.  S\'eminaire sur les Pinceaux de Courbes
Elliptiques (Paris, 1988). Ast\'erisque No. 183 (1990), 7--18.

\spacing \item{[VO1]} Vojta, Paul, Diophantine approximations and
value distribution theory. Lecture Notes in Mathematics, 1239.
Springer-Verlag, Berlin, 1987. x+132 pp.

\spacing \item{[VO2]} Vojta, Paul, On algebraic points on curves.
Compositio Math. 78 (1991), no. 1, 29--36.

\spacing\item{[WI]} Wiles, Andrew Modular elliptic curves and
Fermat's last theorem. Ann. of Math. (2) 141 (1995), no. 3,
443--551.

\spacing \item{[YA1]} Yamanoi, Katsutoshi, On the truncated small
function theorem in Nevanlinna theory. Internat. J. Math. 17 (2006),
no. 4, 417--440.

\spacing \item{[YA2]} Yamanoi, Katsutoshi, Defect relation for
rational functions as targets. With errata by the author. Forum
Math. 17 (2005), no. 2, 169--189.

\spacing \item{[YA3]} Yamanoi, Katsutoshi, The second main theorem
for small functions and related problems. Acta Math. 192 (2004), no.
2, 225--294.

\

\

\spacing

C.~Gasbarri:  Dipartimento di Matematica dell'Universit\`a di Roma
``Tor Vergata", Viale della Ricerca Scientifica, I-00133 Roma (I).

E--mail: gasbarri@mat.uniroma2.it

\bye